\documentclass[10pt]{amsart}
\usepackage[utf8]{inputenc}
\usepackage[T1]{fontenc}
\usepackage{amsfonts,amsmath,amssymb,amsthm,tikz}
\usepackage{graphicx}

\newtheorem{theorem}{Theorem}[section]

\theoremstyle{definition}
\newtheorem{defin}[theorem]{Definition}

\newtheorem{prop}[theorem]{Proposition}
\newtheorem{cor}[theorem]{Corollary}

\newtheorem{thm}[theorem]{Theorem}

\theoremstyle{remark}

\numberwithin{equation}{section}

\begin{document}

\title{Approximate Roots}
\author{Patrick Popescu-Pampu}
\address{Universit\'e de Lille, CNRS, Laboratoire Paul Painlev\'e, 59000 Lille, France}
\email{patrick.popescu-pampu@univ-lille.fr}

This paper appeared in {\em Valuation Theory and its Applications}. Vol. II, F.-V. Kuhlmann, S.Kuhlmann, M.Marshall eds., Fields Institute Communications {\bf 33}, AMS 2003, 285--321.

\vspace{3mm}

\maketitle

\vspace{-5mm}

\begin{abstract}
Given an integral domain $A$, a monic polynomial $P$ of degree $n$ with coefficients in $A$ and a divisor $p$ of $n$, invertible in $A$, there is a unique monic polynomial $Q$ such that the degree of $P-Q^{p}$ is minimal for varying $Q$. This $Q$, whose $p$-th power best approximates $P$, is called the $p$-\textit{th approximate root} of $P$. If $f \in \mathbf{C}[[X]][Y]$ is irreducible, there is a sequence of \textit{characteristic approximate roots} of $f$, whose orders are given by the singularity structure of $f$. This sequence gives important information about this singularity structure. We study its properties in this spirit and we show that most of them hold for the more general concept of \textit{semiroot}. We show then how this local study adapts to give a proof of Abhyankar-Moh's embedding line theorem.
\end{abstract}

\section{Introduction}

The concept of {\it approximate root} was introduced and studied in \cite{AM
  73} in order to prove (in \cite{AM 75}) what is now called the 
  Abhyankar-Moh-Suzuki theorem: it states that the affine line
  can be embedded in a unique way (up to ambient automorphisms) in the affine
  plane. More precisely,
  formulated algebraically the theorem is: 

\medskip

\textbf{Theorem} {\sl (Embedding line theorem)

  If $\mathbf{C}[X,Y] \rightarrow \mathbf{C}[T]$ is an epimorphism of
  $\mathbf{C}$-algebras, then there exists an isomorphism of
  $\mathbf{C}$-algebras $\mathbf{C}[U,V] \rightarrow \mathbf{C}[X,Y]$
  such that the composed epimorphism $\mathbf{C}[U,V] \rightarrow 
  \mathbf{C}[T]$ is given by $U =T, V =0$}.

\medskip

This theorem, as well as other theorems about the group of
automorphisms of $\mathbf{C}[X,Y]$, was seen to be an easy consequence
of the following one, in which $d(P)$ denotes the degree of the
polynomial $P$:

\medskip
    \textbf{Theorem} {\sl (Epimorphism theorem)

 If $\mathbf{C}[X,Y] \rightarrow \mathbf{C}[T]$ is an epimorphism of 
$\mathbf{C}$-algebras, given by $X = P(T), Y = Q(T)$, with 
$d(P)>0, \; d(Q)>0$, 
then $d(P)$
divides $d(Q)$ or vice-versa.}
\medskip

Sometimes in the literature the names of the two theorems are
permuted. The initial proofs (\cite{AM 75}) were simplified in \cite{A 77}. Let
us indicate their common starting point. 

In order to prove the embedding line theorem, Abhyankar and Moh introduced the
image curve of the embedding, whose equation is obtained by computing a
resultant: $f(X,Y)=\mbox{Res}_{T}(P(T)-X, Q(T)-Y)$. The curve $f(X,Y)=0$ has
only \textit{one place at infinity} (see the general algebraic
definition in \cite{A 77}; in our context it means simply that the
closure of the curve in the projective plane has only one point on the
line at infinity and it is unibranch at that point). The fact that
$\mathbf{C}[X,Y] \rightarrow 
\mathbf{C}[T]$ is an epimorphism is equivalent with the existence of a
relation $T=\Psi(P(T),Q(T))$, where $\Psi \in \mathbf{C}[X,Y]$. This
in turn is equivalent with the existence of $\Psi$ such that the
degree of $\Psi(P(T),Q(T))$ is equal to $1$. Now, when $\Psi$ varies,
those degrees form a semigroup. This semigroup was seen to be linked
with a semigroup of the unique branch of $\Psi$ at infinity, which
has a local definition. That is how one passes from a global problem
to a local one.

To describe the situation near the point at infinity, in \cite{AM 75}
the affine plane 
was not seen geometrically as a
chart of the projective plane. The operation was done algebraically,
making the change of variable $X \rightarrow \frac{1}{X}$. So from the
study of the polynomial $f$ one passed to the study of 
$\phi(X,Y) = f(X^{-1},Y)$, seen as an element of
$\mathbf{C}((X))[Y]$. That is why in \cite{AM 73} the local study was
made 
for {\it meromorphic curves}, i.e., elements of
$\mathbf{C}((X))[Y]$. The classical Newton-Puiseux expansions were
generalized to that situation (see the title of \cite{AM 73}) as well
as the notion of semigroup. In order to study this semigroup some 
special approximate roots of $\phi$ were
used, which we call
{\it characteristic approximate roots}. Their importance in this context
lies in the fact that they can be defined globally in the plane, their
local versions being obtained with the same change of variable as
before: $X \rightarrow \frac{1}{X}$.

The proofs in \cite{AM 73} or in \cite{A 77} of the local properties of 
approximate
roots dealt exclusively with locally irreducible
meromorphic curves. 
In \cite{M 82} a generalization for
possibly reducible polynomials was achieved, over an arbitrary non-archimedean
valued field.

An introduction to
Abhyankar's philosophy on curves  and to his
notations can be found in \cite{A 90}. 
A gradual presentation of the general path of the proof of the
epimorphism theorem was tried at
an undergraduate level in \cite{A 75}. See also the presentation done
in \cite{P 77}. Other applications to global
problems in the plane are given in \cite{A 78}. We quote here the
following generalization of the embedding line theorem:
 
\medskip

\textbf{Theorem} {\sl (Finiteness theorem)

Up to isomorphisms of the affine plane, there are only finitely many
embeddings of a complex irreducible algebraic curve with one place at
infinity in the affine plane.}

\medskip

The reference \cite{A 78} also contains some conjectures in
higher dimensions.

Here we discuss mainly the local aspects of approximate roots. We
work in less generality, as suggested by the presentation of the
subject made in \cite{GP 95}. Namely, we consider only
polynomials in $\mathbf{C}[[X]][Y]$. This framework has the
advantage of giving more geometrical insight, many computations being
interpreted in terms of intersection numbers (see also \cite{CW 93}), 
a viewpoint that is
lacking in the meromorphic case. This has also the advantage of
allowing us to interpret the local properties of approximate roots in
terms of the minimal resolution of $f$, a concept which has no analog
in the case of meromorphic curves. We define the concept of
\textit{semiroots}, as being those curves that have the same
intersection-theoretical properties as the characteristic approximate
roots, and we show that almost all the local properties usually used for the
characteristic roots are in fact true for semiroots.

First we introduce the notations for Newton-Puiseux parameterizations of
a plane branch
in arbitrary coordinates, following \cite{GP 95}. 
In section 3 we introduce the general notion of approximate roots. We explain
the concept of semigroup of the branch and related notions in section
4. In section 5  
we introduce the characteristic approximate roots of the branch, we state their
main intersection-theoretical local properties (Theorem \ref{Mainthm})
and we add some  
corollaries. In section 6 we explain the main steps of the proof of
Theorem \ref{Mainthm}. In sections 7 and 8 we give \textit{the proofs} of 
Theorem \ref{Mainthm}, its corollaries and the auxiliary propositions
stated in the text. We prefer to isolate the proofs from the main
text, in order to help reading it. 
In the final section we indicate the changes one must make to the theory
explained before in order to deal with the meromorphic curves and
we sketch a proof of the embedding line theorem.

 A
forerunner of the concept of approximate root was introduced in an
arithmetical context in \cite{ML 360} and \cite{ML 36} (see also
\cite{K 79} for some historical remarks on those papers). The existence
of approximate roots is the content of exercise 13, $\S 1$, in
\cite{B 59}. 
The concept of semiroot is closely associated with that of curve
having maximal contact with the given
branch, introduced in \cite{LJ 69} and \cite{LJ 73}. More details on this
last concept are given in the comments following Corollary
\ref{dual}. Approximate roots of elements of
$\mathbf{C}[[X]][Y]$ are also used in \cite{AO 96} to study the local 
topology of
plane curves. The approximate roots of curves in positive
characteristic are studied in \cite{R 80} using Hamburger-Noether
expansions and the epimorphism theorem is generalized to this case
under some restrictions. The approximate roots of meromorphic curves
are used in 
\cite{A 99} for the study of affine curves with only one irregular
value. The projectivized approximate roots of a curve with one
place at infinity are used in \cite{CPR 00} in order to obtain global
versions of Zariski's theory of complete ideals. In \cite{GP 00}, the
 theorem \ref{Mainthm} proved below and some of its corollaries are
 generalized  
to the case of quasi-ordinary singularities of hypersurfaces. 

We would like to thank S.S.Abhyankar for the explanations he gave us
in Saskatoon on approximate roots. We were also greatly helped in our
learning of the subject by the article \cite{GP 95} of
J.Gwo\'{z}dziewicz and A.P\l oski. We thank
B.Teissier, E.Garc\'{\i}a Barroso and P.D.Gonz\'alez P\'erez for their
comments on preliminary versions of this work and S.Kuhlmann and
F.V.Kuhlmann for the invitation to talk on this subject in Saskatoon.

\section{Notations}

In what follows we do not care about maximal generality on the base
field. We work over $\mathbf{C}$, the field of complex
numbers. By $''a | b''$ we mean ``$a$ divides $b$'', whose negation we
note $''a \not| b''$. The greatest common divisor of
$a_{1},...,a_{m}$ is denoted $\mbox{gcd}(a_{1},...,a_{m})$. If
$q \in \mathbf{R}$, its integral part is denoted $[q]$.

We consider $f(X,Y) \in  \mathbf{C}[[X]][Y]$, a polynomial in
the variable $Y$, monic and irreducible over $\mathbf{C}[[X]]$, the 
ring of formal series in $X$:

$$ f(X,Y) = Y^{N} + \alpha_{1}(X)Y^{N-1} + \alpha_{2}(X)Y^{N-2}+\cdots+
\alpha_{N}(X)$$
where $\alpha_{1}(0)=\cdots=\alpha_{N}(0)=0$.

If we embed $\mathbf{C}[[X]][Y] \hookrightarrow \mathbf{C}[[X,Y]]$, 
the equation $f(X,Y) = 0$ defines a germ of formal (or algebro\"{\i}d)
irreducible curve at the
origin - we
 call it shortly {\it a branch} - in the plane of coordinates
$X,Y$. We denote this curve by $C_{f}$. 

Conversely, if a branch $C \hookrightarrow \mathbf{C}^{2}$ is given, 
the Weierstrass preparation theorem shows that it can be defined by a
unique polynomial of the type just discussed, once the ambient
coordinates $X,Y$ have been chosen, with the exception of
$C=Y$-axis. If we describe like this a curve 
by a polynomial equation $f$ with respect to the variable $Y$, we call
briefly its degree $N$ in $Y$ {\it the degree of} $f$, and we denote it by 
$d(f)$ or $d_{Y}(f)$ if we want to emphasize the variable in which it
is polynomial. When $C$ is {\it transverse} to the $Y$-axis (which means that 
the 
tangent cones of $C$ and of the $Y$-axis have no common components), 
we have the equality $d(f)=m(C)$, where
$m(C)$ denotes the multiplicity of $C$ (see section \ref{Semigroup}).

From now on, each time we speak about the curve $C$, we mean the curve
$C_{f}$, for the fixed $f$.

The curve $C$ can always be parameterized in the following way (see
\cite{EC 15}, \cite{Z 35}, 
\cite{A 90}, \cite{T 95}, \cite{W 50}):
\begin{equation}
\label{Puiseux}
 \begin{cases}
   X = T^{N} \\
   Y = \sum_{j \geq 1} a_{j}T^{j} = \cdots+ a_{B_{1}}T^{B_{1}} + \cdots +
       a_{B_{2}}T^{B_{2}} + \cdots + a_{B_{G}}T^{B_{G}} + \cdots
 \end{cases}
\end{equation}
with $\mbox{gcd}(\{N\} \cup \{j,a_{j} \neq 0 \}) = 1$.

The exponents $B_{j}$, for $j\in \{ 1,...,G\}$ are defined inductively:

\hspace{5mm}
 $B_{1}: = \mbox{min} \{j, a_{j} \neq 0, N \not| j\}$

\hspace{5mm}
 $B_{i}: = \mbox{min} \{j, a_{j} \neq 0, \mbox{gcd}(N,B_{1},...,B_{i-1})
 \not| j\}$, for $i \geq 2. $

The number $G$ is the least one for which $\mbox{gcd}(N,B_{1},...,B_{G})=1$.

We define also: $B_{0} := N =d(f)$. Then $(B_{0},B_{1},...,B_{G})$ is called 
the {\it characteristic sequence} of $C_{f}$ in the
coordinates $X,Y$. The $B_{i}$'s are {\it the characteristic
  exponents} of $C_{f}$ {\it with respect to} $(X,Y)$.

 A parameterization like (\ref{Puiseux}) is called {\it a primitive 
Newton-Puiseux  parameterization with respect to} $(X,Y)$ of
the plane branch $C$. Notice that $X$ and $Y$ cannot be permuted in
  this definition. 

Let us explain why we added the attribute
``primitive''. If we write $T=U^{M}$, where $M \in \mathbf{N}^*$, we
obtain a parameterization using the variable $U$. In the new
parameterization, the greatest common divisor of the exponents of the
series $X(U)$ and $Y(U)$ is no longer equal to $1$. In this case we
say that the
parameterization is {\it not primitive}.  When we speak only of a 
``Newton-Puiseux parameterization'',
we mean a primitive one.

We define now the sequence of greatest common divisors:
$(E_{0}, E_{1},...,E_{G})$ in the following way:
$$  E_{j} = \mbox{gcd}(B_{0},...,B_{j}) \mbox{ for } j \in \{0,...,G\}.$$ 
In particular:
$  E_{0} = N, \: E_{G}=1$. 
Define also their quotients:
$$ N_{i} = \frac{E_{i-1}}{E_{i}}>1, \makebox{ for } 1 \leq i \leq G. $$
This implies: 
$$E_{i} = N_{i+1}N_{i+2}\cdots N_{G}, \makebox{ for } 0 \leq i \leq
G-1.$$

\vspace{5mm}

Let us introduce the notion of {\it Newton-Puiseux series of}  $C$ 
{\it with respect to} $(X,Y)$. It is a series of the form:

\begin{equation}
  \label{Puiseries}
 \eta(X) = \sum_{j \geq 1} a_{j}X^{\frac{j}{N}} 
\end{equation}

\noindent
obtained from (\ref{Puiseux}) by replacing $T$ by $X^{\frac{1}{N}}$. It is
an element of $\mathbf{C}[[X^{\frac{1}{N}}]]$. One has then the
equality $f(X,\eta(X))=0$, so $\eta(X)$ can be seen as an expression for {a
root} of the polynomial equation in $Y$: $f(X,Y)=0$. All the other
roots can be obtained from (\ref{Puiseries}) by changing  $X^{\frac{1}{N}}$ 
to $\omega X^{\frac{1}{N}}$, where  $\omega$ is an arbitrary $N$-th
root of unity. This is a manifestation of the fact that the Galois
group of the field extension 
$\mathbf{C}((X))\rightarrow \mathbf{C}((X^{\frac{1}{N}}))$ is 
$\mathbf{Z}/N\mathbf{Z}$. From this remark we get another
presentation of characteristic exponents: 

\begin{prop}
\label{differences}
 The set $\{ \frac{B_{1}}{N},...,\frac{B_{G}}{N}\}$ is equal to the
 set: $$\{v_{X}(\eta(X)-\zeta(X)), \mbox{ } \eta(X) \mbox{ and }
 \zeta(X) \mbox{ are distinct roots of } f \}.$$
\end{prop}

\medskip

Here $v_{X}$ designates the order of a formal fractional power series in the
variable $X$.

Given a Newton-Puiseux series (\ref{Puiseries}), define for $k \in
\{0,...,G\}$:

\[ \eta_{k}(X) = \sum_{1 \leq j < B_{k+1}} a_{j}X^{\frac{j}{N}}. \]

It is the sum of the terms of $\eta(X)$ of exponents strictly less than 
$\frac{B_{k+1}}{N}$. We call $\eta_{k}(X)$ {\it a k-truncated Newton-Puiseux
  series of} $C$ {\it with respect to} $(X,Y)$. If the
 parameterization (\ref{Puiseux}) is reduced, then $\eta_{k}(X) \in 
 \mathbf{C}[[X^{\frac{E_{k}}{N}}]]$ and there are exactly
$\frac{N}{E_{k}}$ such series.
 
\vspace{5mm}
Before introducing the concept of approximate root, we give an
example of a natural question about Newton-Puiseux parameterizations, which
will be answered very easily using that concept.

\vspace{5mm}

\noindent \textbf{Motivating example}

There are algorithms to compute Newton-Puiseux parameterizations of the 
branch 
starting from the polynomial $f$. If one wants to know only the
beginning of the parameterization, one could ask if it is enough to know
only some of the coefficients of the polynomial $f$. The answer is
affirmative, as is shown by the following proposition:

\begin{prop}
\label{example} If $f$ is irreducible with characteristic sequence 
$(B_{0},...,B_{G})$, then 
the terms of the $k$-truncated Newton-Puiseux series of $f$ 
depend only on $\alpha_{1}(X),...,\alpha_{\frac{N}{E_{k}}}(X)$.
\end{prop}

\medskip

The proof will appear to be very natural once we know the concept of
approximate root and Theorem \ref{Mainthm}. Let us illustrate the
proposition by a concrete case.

Consider:

\hspace{5mm}
$f(X,Y) = Y^{4}-2X^{3}Y^{2}-4X^{5}Y+X^{6}-X^{7}$.

One of its Newton-Puiseux parameterizations is:

\hspace{5mm}
$ \begin{cases}
    X = T^{4} \\
    Y = T^{6} + T^{7}
  \end{cases} $

We get, using the proposition for $k=1$, that 
every irreducible polynomial of the form:
$$g(X,Y) = Y^{4}-2X^{3}Y^{2} + \alpha_{3}(X)Y + \alpha_{4}(X)$$
whose characteristic sequence is $(4,6,7)$, has a Newton-Puiseux 
series of the type:
$$ Y = X^{\frac{3}{2}}+ \gamma(X),$$
with $v_{X}(\gamma) \geq \frac{7}{4}$. 

\medskip 

It is now the time to introduce the approximate roots...

\section{The Definition of Approximate Roots} \label{Defroots}

Let $A$ be an integral domain (a unitary commutative ring without zero
divisors). If $P \in A[Y]$ is a polynomial
with coefficients in $A$, we shall denote by $d(P)$ its degree. 

Let $P \in A[Y]$ be monic of degree $d(P)$, and $p$ a divisor of 
$d(P)$.
In general there is no polynomial $Q \in A[Y]$ such that $P = Q
^{p}$, i.e. there is no {\it exact root} of order $p$ of the
polynomial $P$. But one can ask for a {\it best approximation} of this
equality. We speak here of approximation in the sense that the
difference $P -Q^{p}$ is of degree as small as possible for varying
$Q$. Such a $Q$ does not necessarily exist. But it exists if one has
the following
condition on the ring $A$, verified for example in the case that
interests us here, $A = \mathbf{C}[[X]]$:
$p$ {\it is invertible in} $A$.

More precisely, one has the following proposition:

\begin{prop} \label{defrt}
 If $p$ is invertible in $A$ and $p$ divides $d(P)$, then there is a unique
 monic polynomial $Q \in A[Y]$ such that: 
 \begin{equation}
   \label{ineg} 
   d(P -Q^{p}) < d(P)-\frac{d(P)}{p}.
 \end{equation}
\end{prop}

This allows us to define:

\begin{defin} \label{approx}
 The unique polynomial $Q$ of the preceding proposition is named the 
$\mathbf{p}$-\textbf{th approximate root} of $P$. It is denoted 
$\mathbf{\sqrt[p]{P}}$.
\end{defin}

Obviously:
$$ d(\sqrt[p]{P})=\frac{d(P)}{p}.$$
\noindent \textbf{Example}: Let $P = Y^{n} + \alpha_{1} Y^{n-1} + \cdots
+ \alpha_{n}$ 
be an element of $A[Y]$. Then, if $n$ is invertible in $A$:
$$ \sqrt[n]{P} = Y + \frac{\alpha_{1}}{n}.$$
We recognize here the {\it Tschirnhausen transformation} of the
variable $Y$. That is the reason why initially (see the title of
\cite{AM 73}) the approximate roots were seen as generalizations of
the Tschirnhausen transformation.

We give now a proposition showing that in some sense the
notation $\sqrt[p]{P}$ is adapted:

\begin{prop}
  \label{tworoots}
  If $p,q \in \mathbf{N}^*$ are invertible in $A$, then $\sqrt[q]{\sqrt[p]{P}}
  =\sqrt[pq]{P}$.
\end{prop}

We see that approximate roots behave in this respect like usual $d$-th
roots. The following construction shows another link between the two
notions. We add it for completeness, since it will not be used in the
sequel.

Let $P \in A[Y]$ be a monic polynomial. Consider $P_{1} \in
A[Z^{-1}]$, $P_{1}(Z)= P(Z^{-1})$. If we embed the ring $A[Z^{-1}]$
into $A((Z))$, the ring of meromorphic series with coefficients in
$A$, the $p$-th root of $P_{1}$ exists inside $A((Z))$. It is the
unique series $P_{2}$ with principal term $1 \cdot Z^{-\frac{n}{p}}$ 
such that $P_{2}^{p} = P_{1}$. We note:
$$ P_{1}^{\frac{1}{p}}:= P_{2}. $$

Consider the {\it purely meromorphic part} $M(P_{1}^{\frac{1}{p}})$ of 
$P_{1}^{\frac{1}{p}}$, the sum of the terms having $Z$-exponents $\leq
0$.

We have $M(P_{1}^{\frac{1}{p}}) \in A[Z^{-1}] $, so: 
$$Q(Y)= M(P_{1}^{\frac{1}{p}})(Y^{-1}) \in A[Y].$$ We can state now the 
proposition (see \cite{M 80}, \cite{M 82}):

\begin{prop}
   \label{second}
If $Q \in A[Y]$ is defined as before, then $Q = \sqrt[p]{P}$.
\end{prop}

\section{The Semigroup of a Branch} \label{Semigroup}

Let $\mathcal{O}_{C}= \mathbf{C}[[X]][Y]/(f)$ be the local ring of the
germ $C$ at the origin. It is an integral local ring of
dimension 1. 

Let $\mathcal{O}_{C} \rightarrow
\overline{\mathcal{O}}_{C}$ be the morphism of normalization of 
$\mathcal{O}_{C}$, i.e., $\overline{\mathcal{O}}_{C}$ is the integral
closure of $\mathcal{O}_{C}$ in its field of fractions. This new  ring 
is regular (normalization is a desingularization in dimension 1), and
so it is {\it a discrete valuation ring of rank} 1. Moreover,
there exists an element $T \in \overline{\mathcal{O}}_{C}$ , called
{\it a
uniformizing parameter}, such that  $\overline{\mathcal{O}}_{C} \simeq 
\mathbf{C}[[T]]$. Then the valuation is simply the $T$-adic valuation 
$v_{T}$, which associates to each element of
$\overline{\mathcal{O}}_{C}$, seen as a series in $T$, its order in $T$. 

\begin{defin}
  The \textbf{semigroup} $\Gamma(C)$ of the branch $C$ is the image by the
  $T$-adic valuation of the non zero elements of the ring $\mathcal{O}_{C}$:
 \[  \Gamma(C) := v_{T}(\mathcal{O}_{C}-\{0\}) \subset 
v_{T}(\overline{\mathcal{O}}_{C}-\{0\}) = \mathbf{N} = \{0,1,2,...\}. \]
\end{defin}

\medskip

The set $\Gamma(C)$ is indeed a semigroup, which comes from the
additivity property of the valuation $v_{T}$:
\[ \forall \phi,\psi \in \mathcal{O}_{C}-\{0\}, v_{T}(\phi
\psi)=v_{T}(\phi)+v_{T}(\psi) .\]

The previous definition is intrinsic, it does not depend on the fact
that the curve $C$ is planar. Let us now turn to other interpretations
of the semigroup. 

First, our curve is given with a fixed embedding in the plane of
coordinates $(X,Y)$. Once we have chosen a uniformizing parameter $T$,
we have obtained a {\it parameterization} of the curve: 
$ \left\{ \begin{array}{l} X=X(T) \\ Y=Y(T)
         \end{array} 
  \right.$.
   
For example, a Newton-Puiseux parameterization would work.

If $f' \in \mathcal{O}_{C}-\{0\}$, it can be seen as the restriction of an
element of the ring $\mathbf{C}[[X]][Y]$, which we denote by
the same symbol $f'$. The curve $C'$ defined by the equation $f'=0$ has
an intersection number with $C$ at the origin. We note it 
$(f,f')$, or $(C,C')$, to insist on the fact that this
number depends only on the curves, and not on the coordinates or the
defining equations. 
We have then the equalities:
\[ v_{T}(f') = v_{T}(f'(X(T),Y(T))) = (f,f') ,\]
which provides a geometrical interpretation of the semigroup of
the branch $C$:
\[ \Gamma(C) = \{ (f,f'), f' \in \mathbf{C}[[X]][Y], f \not\vert f' \}. \]
From this viewpoint, the semigroup is simply the set of possible
intersection numbers with curves not containing the given branch. 

The minimal non-zero element of $\Gamma(C)$ is the \textit{multiplicity}
$m(C)$, noted also $m(f)$ if $C$ is defined by $f$. It is the lowest
degree of a monomial appearing in the Taylor series of $f$, and
therefore also the 
intersection number of $C$ with smooth curves passing through the
origin and transverse to the tangent cone of $C$.

If $p_{1},...,p_{l}$ are elements of $\Gamma(C)$, the sub-semigroup
 $\mathbf{N}p_{1}+\cdots +\mathbf{N}p_{l}$ they generate is denoted by: 
\[ \langle p_{1},...,p_{l} \rangle. \]
One has then the following result,
expressing a set of generators of the semigroup in terms of the
characteristic exponents:

\begin{prop}
 \label{semigroup}
 The degree $N$ of the polynomial $f$ is an element of $\Gamma(C)$, 
denoted $\overline{B}_{0}$. So, $\overline{B}_{0}=B_{0}$. Define
inductively other numbers $\overline{B}_{i}$ 
by the following property:
 \[ \overline{B}_{i} = \mbox{min}\{j \in \Gamma(C), j \notin 
  \langle \overline{B}_{0},...,\overline{B}_{i-1} \rangle \} .\]
Then this sequence has exactly $G+1$ terms 
$\overline{B}_{0},...,\overline{B}_{G}$, which verify the following
properties for $0 \leq i \leq G$ (we consider by definition that 
$\overline{B}_{G+1} = \infty$):

\[ \begin{array}{l}
      1) \: \overline{B}_{i} = B_{i} + 
            \sum_{k=1}^{i-1}\frac{E_{k-1}-E_{k}}{E_{i-1}} B_{k}. \\
      2) \: \mbox{gcd}(\overline{B}_{0},...,\overline{B}_{i}) = E_{i}. \\
      3) \:  N_{i} \overline{B}_{i} < \overline{B}_{i+1}.
   \end{array}
\]
\end{prop}
\medskip 

A proof of this proposition is given in \cite{Z 73} for generic
coordinates (see the definition below) and in \cite{GP 95} in this
general setting. Other properties
of the generators are given in \cite{M 77}, in the generic case (see
the definition below). In fact the proof can be better conceptualized
if one uses the notion of semiroot, more general than the notion of
characteristic root. This notion is explained in section 6. When
reading the proof of Proposition  \ref{semigroup}, one should become convinced
that there is no vicious circle in the use of the $\overline{B}_{k}$'s.

The point is that it appears easier to define the $\overline{B}_{k}$'s
by property 1) and then to prove the minimality property of the
sequence. We have used the other way round in the formulation of
Proposition  \ref{semigroup} because at first sight the minimality definition
seems more natural. One should also read the comments preceding 
Proposition  \ref{strict}. 

 The generators of the semigroup introduced in this
proposition depend on the coordinates $X,Y$, but only in a loose
way. Indeed, they are uniquely determined by the semigroup once the
first generator $\overline{B}_{0}$ is known. This generator, being
equal to the degree of the polynomial $f$, depends on
$X,Y$. Geometrically, 
it is the intersection number $(f,X)$. It follows from this that for
{\it generic coordinates}, i.e. with the $Y$-axis transverse to the
curve $C$, the generators are independent of the coordinates. 

We can therefore speak of {\it generic characteristic exponents}. They are
a complete set of invariants for the equisingularity and the
topological type of the branch (see \cite{Z 73}). For the other
discrete invariants introduced before, we speak in the same way of
\textit{generic} ones and we use lower case letters to denote them, as
opposed to capital ones for the invariants in arbitrary
coordinates. Namely, we use the following notations for the
generic invariants:
\[  \begin{array}{l}
        n\\
 (b_{0},...,b_{g})\\
  (e_{0},...,e_{g})\\
  (n_{1},...,n_{g})\\
  (\overline{b}_{0},...,\overline{b}_{g})
     \end{array}.\]
This makes it easy to recognize in every context if we suppose the
coordinates to be generic or not.

We call $g$ the 
\textit{genus} of the curve $C$. 

The exponent $b_{0}$ is equal to the multiplicity $m(C)$ of
$C$ at the origin. When $D$ is a curve passing through $0$, we have 
$(C,D)=b_{0}$ if and only if $D$ is smooth and transverse to $C$
at $0$. 

The
preceding proposition shows that in arbitrary coordinates the
generators of the semigroup  are determined by the
characteristic exponents. But the relations can be reversed, and show that
conversely, the characteristic exponents are determined by the generators of
the semigroup. From this follows the invariance of the characteristic 
exponents with respect to the generic coordinates chosen for the
computations. Moreover, like this one can easily obtain a proof of the
classical {\it inversion formulae} for plane branches (see another
proof in \cite{A 67}). Let us state it in a little extended form. 

Let $(X,Y)$ and $(x,y)$ be two systems of coordinates, the
second one being generic for $C$. We consider the characteristic
exponents $(B_{0},...,B_{G})$ of $C$ with respect to $(X,Y)$.

\begin{prop}
  \label{Inversion}
 (Inversion formulae) 

The first characteristic exponent $B_{0}$ can take values only in the
set \\$\{ lb_{0}, 1 \leq l \leq [\frac{b_{1}}{b_{0}}]\} 
\cup \{ b_{1}\}$. The knowledge of its value completely determines
the rest of the exponents in terms of the generic ones: 

1) $B_{0} = b_{0} \Rightarrow G=g$ and: $$(B_{0},...,B_{G})=
(b_{0},...,b_{g}).$$

2) $B_{0} =  lb_{0}$ with $2 \leq l \leq
[\frac{b_{1}}{b_{0}}] 
\Rightarrow G = g+1$ and: 
$$(B_{0},...,B_{G})= (lb_{0}, b_{0}, 
b_{1}+(1-l)b_{0},..., b_{g}+(1-l)b_{0}).$$ 

3) $B_{0} =b_{1} \Rightarrow G=g$ and: 
$$ (B_{0},...,B_{G})= 
(b_{1}, b_{0}, b_{2}+b_{0}-b_{1},
b_{3}+b_{0}-b_{1},...,b_{g}+b_{0}-b_{1}).$$ 

Moreover, for $k \in \{1,...,g \}$, the $k$-truncations of the Newton-Puiseux
series with respect to $(x,y)$ depend only on the
$(k+\epsilon)$-truncation of the Newton-Puiseux series with respect to
$(X,Y)$, where $\epsilon=G-g \in \{0,1\}$.

\end{prop}

The name given classically to one form or another of this proposition
comes from the fact it answers the question: {\it what can we say
  about the Newton-Puiseux series with respect to} $(Y,X)$ {\it if we know
  it in terms of} $(X,Y)$? In this question, one simply {\it inverts}
  the coordinates. 

We prove the statement on truncations using, as in the case of
Proposition  \ref{semigroup}, the notion of \textit{semiroot},
introduced in section 6.

\section{The Main Theorem}

As before, the polynomial $f \in \mathbf{C}[[X]][Y]$ is supposed to be
irreducible. To be concise, we note in what follows:
 
$$f_{k} := \sqrt[E_{k}]{f}.$$

Next theorem is the main one, (7.1), in \cite{AM 73}. A different
proof is given in \cite{A 77}, (8.2). Here we give a proof inspired by
\cite{GP 95}.

\begin{thm}
 \label{Mainthm}
The approximate roots $f_{k}$ for $0\leq k \leq G$, have the following
properties:

\hspace{5mm} 1) $d(f_{k}) = \frac{N}{E_{k}}$ and 
                $(f, f_{k}) = \overline{B}_{k+1}$.

\hspace{5mm} 2) The polynomial $f_{k}$ is irreducible and its
characteristic exponents in these coordinates are $\frac{B_{0}}{E_{k}}, 
\frac{B_{1}}{E_{k}},...,\frac{B_{k}}{E_{k}}$.
\end{thm}

Theorem \ref{Mainthm} gives properties of some of the approximate roots of
$f$. One does not consider all the divisors of $N$, but only some
special ones, computed from the knowledge of the characteristic 
exponents. For this reason, we name them {\it
  the characteristic approximate roots of} $f$. 

We give now a list of corollaries. In fact these corollaries hold 
more generally for \textit{semiroots}, see the comments made after Definition
\ref{semiroot}. The proofs of the theorem and of its corollaries are
given in section \ref{Proofs}. Before that, in section \ref{Steproof}
we explain the main steps in the proof of Theorem \ref{Mainthm}
 
\begin{cor}
\label{comput}
 The irreducible polynomial $f$ being given, one can compute
 recursively its characteristic approximate roots in the following
 way. Compute the $N$-th root $f_{0}$ of $f$ and put $E_{0}=N$. 
If $f_{k}$ was 
computed, put  $(f, f_{k}) = \overline{B}_{k+1}$. As 
$E_{k}$ has already been computed,
 take $E_{k+1} = \mbox{gcd}(E_{k},\overline{B}_{k+1})$ and compute 
$f_{k+1} = \linebreak =
\sqrt[E_{k+1}]{f}$. One can then deduce the characteristic 
 exponents from the characteristic roots.
\end{cor}

 This has been extended to the case of meromorphic curves in \cite{A
   94}. The preceding algorithm
works only if $f$ is irreducible. But it can be adapted to give a
method of deciding whether a given $f$ is indeed irreducible, as was done
in \cite{A 89}. See also the more elementary presentation given in
\cite{A 88}. A generalization of this criterion of irreducibility to
the case of arbitrary characteristic is contained in \cite{CM 00}.

Following the proof of the proposition, we add an example of
application of the algorithm.
 
\begin{cor}
 \label{coincid}
For $0 \leq k \leq G$, the polynomials $f$ and $f_{k}$ have equal
sets of $k$-truncations of their
Newton-Puiseux series.
\end{cor}

This, together with the remark following equation (\ref{recursion}),
gives an immediate proof of Proposition \ref{example}.

\begin{cor}
 \label{expansions}
 Every $\phi \in \mathbf{C}[[X]][Y]$ can be {\it uniquely} written as
 a finite sum of the
 form:
\[ \phi = \sum_{i_{0},...,i_{G}}\alpha_{i_{0}...i_{G}}f_{0}^{i_{0}} 
f_{1}^{i_{1}}\cdots f_{G}^{i_{G}} \]
where $i_{G} \in \mathbf{N}$,  $0 \leq i_{k} < N_
{k+1}$ for $0 \leq k \leq
G-1$ and the coefficients $\alpha_{i_{0}...i_{G}}$ are elements of the
ring $\mathbf{C}[[X]]$. Moreover:

\hspace{5mm} 1) the $Y$-degrees of the terms appearing in the right-hand
side of the preceding equality are all distinct.

\hspace{5mm} 2) the orders in $T$ of the terms  
$$\alpha_{i_{0}...i_{G}}(T^{N})f_{0}(T^{N},Y(T))^{i_{0}} 
\cdots f_{G-1}(T^{N},Y(T))^{i_{G-1}}$$ 
are all distinct, where $T \rightarrow (T^{N}, Y(T))$ is a Newton-Puiseux
parameterization of $f$.
\end{cor}

 There is no a priori bound on $i_{G}$: this exponent is equal to
 $[\frac{d(\phi)}{N}]$. The orders in $T$ appearing in 2) are the
 intersection numbers of $f$ with the curves 
defined by the terms of the sum which are not divisible by $f$.

This corollary is essential for the applications of Theorem
\ref{Mainthm} to the proof of the
embedding line theorem. Indeed, the point 2) allows one to compute
$(f,\phi)$ in terms of the numbers $(f, \alpha_{i_{0}...i_{G}}f_{0}^{i_{0}} 
f_{1}^{i_{1}}\cdots f_{G}^{i_{G}})$. But, as explained in the
introduction, one is interested in the semigroup of $f$, composed of
the intersection numbers $(f,\phi)$ for varying $\phi$. This way of
studying the semigroup of $f$ is the one focused on in \cite{AM 73}
and \cite{A 77}.

\begin{cor}
 \label{graduation}
 The images of $X,
f_{0},f_{1},...,f_{G-1}$ into the graded ring 
$gr_{v_{T}}\mathcal{O}_{C}$ generate it as a $\mathbf{C}$-algebra.  
If the coordinates are
generic, they form a minimal system of generators.
\end{cor}

We have defined generic coordinates in the remark following 
Proposition \ref{semigroup}.
Here $gr_{v_{T}}\mathcal{O}_{C}$ is {\it the graded ring of}
$\mathcal{O}_{C}$ {\it with respect to the valuation} $v_{T}$. This
concept is defined
in general,
if $A$ is a domain of integrity, $F(A)$ its field of fractions and
$\nu$ a valuation of  $F(A)$ that is positive on $A$. In this
situation, we define first the semigroup of values $\Gamma(A)$ to be 
the image of $A-\{0\}$ by the valuation. If $p \in \Gamma(A)$, we define the
following ideals:
 \[ I_{p} = \{ x \in A, \nu(x) \geq p\}, \]
 \[ I_{p}^{+} = \{ x \in A, \nu(x) > p\}. \]
The graded ring of $A$ with respect to the valuation $\nu$ is defined 
in the following way:
\[ gr_{\nu}A = \bigoplus_{p \in \Gamma(A)}I_{p}/I_{p}^{+}. \]

This viewpoint on the approximate roots is focused on in \cite{S 90} and 
\cite{S 97}, where the general concept of {\it
  generating sequence} for a valuation is introduced. This concept
generalizes the sequence of characteristic approximate roots,
introduced before. 

In the case of irreducible germs of plane curves, the spectrum of
$gr_{v_{T}}(\mathcal{O}_{C})$ is the so-called \textit{monomial curve}
associated to $C$. It was used in \cite{GT 00} in order to show that one
could understand better the desingularization of $C$ by embedding it
in a space of higher dimension.  
\medskip

Before stating the next corollary, let us introduce some other notions.
For more details one can consult \cite{BK 86}, \cite{LMW
  89} and \cite{S 90}.

An {\it embedded resolution} of $C$ is a proper birational morphism 
$ \pi : \Sigma \rightarrow \mathbf{C}^{2}$ such that $\Sigma$ is
smooth and the total transform $\pi^{-1}(C)$ is a divisor with normal
crossings. Such morphisms exist and they all factorize through a
minimal one $\pi_{m} :\Sigma_{m} \rightarrow \mathbf{C}^{2}$ which can
be obtained in the following way. Start from $C \hookrightarrow 
\mathbf{C}^{2}$ and blow-up the origin. Take the total transform
divisor of 
$C$ in the resulting surface. All its points are smooth or with
normal crossings, with the possible exception of the point on the
strict transform of $C$. If at this point the divisor is not with
normal crossing, 
blow up the point. Then repeat the process. After a finite number of steps, one
obtains the minimal embedded resolution of $C$.

\begin{figure}[tb]
\vspace{5pc}
\begin{center}
\includegraphics[width=3.5in]{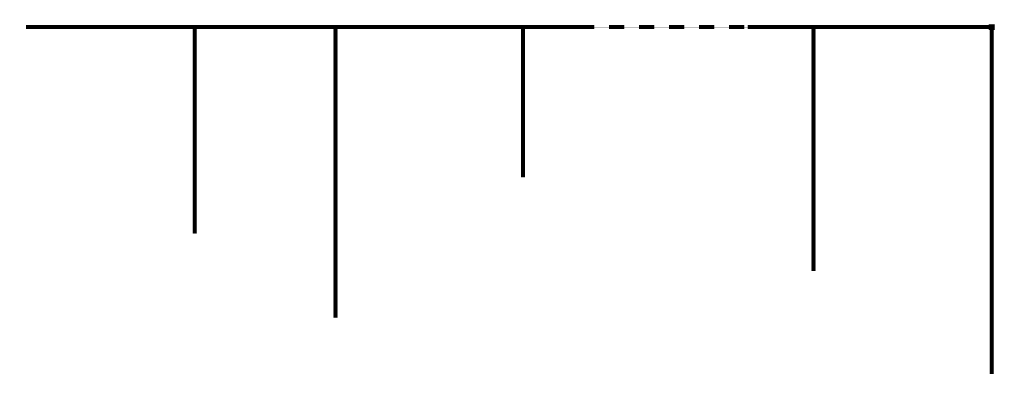}
\end{center}
\caption{The Dual Graph}
\label{dualgr}
\end{figure}

The reduced exceptional divisor $\mathcal{E}$ of $\pi_{m}$ is
connected, which can be easily seen from the previous description by
successive blowing-up. This phenomenon is much more general, and known
under the 
name ``Zariski's connectedness theorem'' or ``Zariski's main
theorem'', see \cite{Z 55}, \cite{M 76} and \cite{H 77}. The
components of $\mathcal{E}$ are isomorphic to 
$\mathbf{C}\mathbf{P}^{1}$. We consider 
{\it the dual graph} $D(\pi_{m})$ of $\mathcal{E}$, whose vertices are
in bijection 
with the components of $\mathcal{E}$. Two vertices are connected by an edge if
and only if the corresponding components intersect on
$\Sigma_{m}$. The graph $D(\pi_{m})$ is then a tree like in Figure
\ref{dualgr}, in which we represent only the underlying topological
space of the graph and not its simplicial decomposition.

In this picture there are exactly $g$ vertical segments, $g$ being the
genus of $f$ (see its definition in the comments following
Proposition \ref{semigroup}). The first vertex on the left of the horizontal 
segment corresponds to the component of $\mathcal{E}$ created by the first 
blowing-up. The vertex of attachment of the horizontal segment and of the 
right-hand vertical segment corresponds to the component of $\mathcal{E}$ 
which cuts the strict transform of $C$.

If we consider also the strict transform of $C$ on $\Sigma_{m}$,
we represent it by an arrow-head vertex connected to the vertex of
$D(\pi_{m})$ which represents the unique
component of $\mathcal{E}$ which it intersects. We denote this new
graph by $D(\pi_{m},f)$.

This graph as well as various numerical characters of the components
of $\mathcal{E}$ can be computed from a generic Newton-Puiseux series for
$f$. The first to have linked Newton-Puiseux series with the resolution of
the singularity seems to be M.Noether in \cite{N 90}. See also
\cite{EC 15} for the viewpoint of the italian school.

\begin{figure}[tb]
\begin{center}
\includegraphics[width=3.5in]{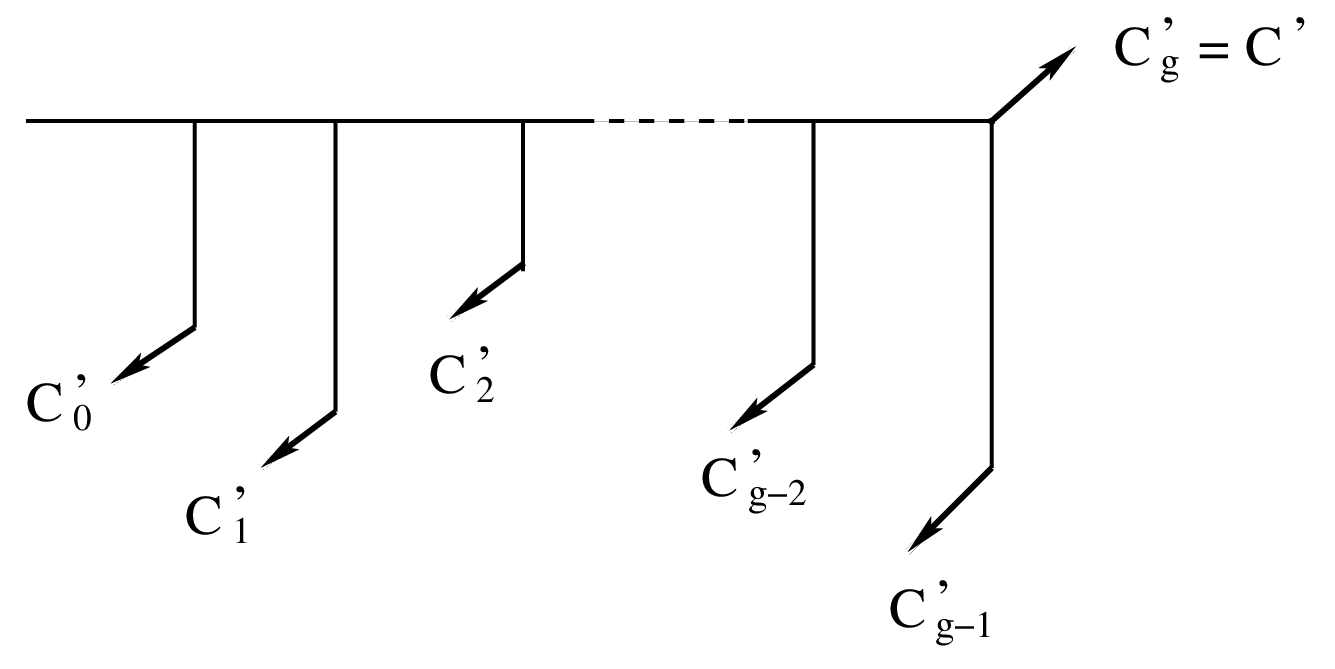}
\end{center}
\caption{The Total Dual Graph}
\label{totaldualgr}
\end{figure}

\begin{cor} \label{dual}
 Let $\pi_{m}$ be  the minimal embedded resolution of
 $C_{f}$. We consider the characteristic approximate roots $f_{k}$,
 for  
$0 \leq k \leq g$ with respect to generic coordinates. Let us denote 
by $C_{k}$ the curve defined by the equation $f_{k}=0$. One has
 evidently $C_{f}= 
C_{g}$. Let us also denote by $C_{k}'$ the strict transform of $C_{k}$
 by the morphism 
$\pi_{m}$. Then the curves $C_{k}'$ are smooth and transverse to a unique
 component of the exceptional divisor of $\pi_{m}$. The dual graph of
 the total transform of $f_{0}f_{1}\cdots f_{g}$ is represented in
 Figure \ref{totaldualgr}.
\end{cor}

The previous corollary gives a topological interpretation of the
characteristic approximate roots, showing how they can be seen as
generalizations 
of smooth curves having maximal contact with $C$. 

Such a generalization was already made in \cite{LJ 69} and \cite{LJ
 73}, where the
notion of {\it maximal contact} with $f$ was extended from smooth
curves to singular curves having at most as many generic
characteristic exponents as $f$. It was further studied in \cite{C
 80}. Let us explain this notion.

If $D$ is a plane branch, let
$$\nu(D):=\frac{1}{m(D)}\mbox{sup}_{D'}\{(D,D')\},$$ where the
supremum is taken over all the choices of smooth $D'$. It is a finite
rational number, with the exception of the case when $D$ is smooth,
which implies $\nu(D)=+ \infty.$

Consider now the sequence of point blowing-ups
which desingularizes $C$. For $i \in \{0,...,g\}$, let $D_{i}$ be the
first strict transform of $C$ that has genus
$g-i$. One has $D_{0}=C$. Define:
$$ \nu_{i}(C):=\nu(D_{i}).$$
The sequence $(\nu_{0}(C),...,\nu_{g}(C))$ was named in \cite{LJ
  69} {\it the sequence of Newton
coefficients} of $C$. In characteristic $0$ - for example when working
over $\mathbf{C}$, as we do in this article - its knowledge is
equivalent to the knowledge of the characteristic sequence. The advantage of
the Newton coefficients is that they are defined in any characteristic.

\begin{defin}
If $D$ is a branch of genus $k \in \{0,...,g\}$, we say that $D$
\textbf{has maximal contact} with $C$ if $ \nu_{i}(D)=\nu_{i}(C)$
for every $i \in \{0,...,k\}$ and $(C,D)$ is the supremum  of the
intersection numbers of $C$ with curves of genus $k$ having the
previous property.
\end{defin}

It can be shown with the same kind of arguments as those used to prove
Corollary \ref{dual}, that for every 
$k \in \{0,...,g\}$, the curves having genus $k$ 
and maximal contact with $f$ are exactly the $k$-semiroots in generic 
coordinates.

In order to understand better Corollary \ref{dual}, let us introduce
another concept:

\begin{defin}
Let $L$ be some component of the reduced exceptional divisor
$\mathcal{E}$. A branch $D \hookrightarrow  
\mathbf{C}^{2}$ is called \textbf{a curvette} with respect to $L$ if its
strict transform by $ \pi_{m}$ is smooth and transversal to $L$ at a
smooth point of $\mathcal{E}$.
\end{defin}

Let $L_{0}$ be the component of $\mathcal{E}$ created by the
blowing-up of $0 \in
\mathbf{C}^{2}$. For every $k \in \{1,...,g\}$, let $L_{k}$ be the
  component at the free end of the $k$-th vertical segment of
  $D(\pi_{m})$. Let $L_{g+1}$ be the component intersecting the strict
  transform of $C$.

\begin{cor}  \label{curvette}
A characteristic approximate root of $f$ in arbitrary coordinates is a 
curvette with
respect to one of the components $L_{0},L_{1},...,L_{g+1}$.
\end{cor}

\medskip

This corollary is an improvement of Corollary \ref{dual}, which says 
this is true in generic coordinates. This more general property is
important for the geometrical interpretations of approximate roots
given in \cite{CPR 00}. A deeper study of curvettes, for possibly
multi-branch curve singularities can be found in \cite{LMW 89}.

\section{The Steps of the Proof} \label{Steproof}

In this section we explain only the main steps in the proof of Theorem 
\ref{Mainthm}, as well as a reformulation for the corollaries. The
complete proofs are given in section \ref{Proofs}.

First we have to introduce a new notion, fundamental for the proof,
that of the {\it expansion of a polynomial} in terms of another
polynomial. This is the notion mentioned in the title of 
\cite{A  77}.

Let $A$ be an integral domain and let $P, Q \in A[Y]$ be monic
polynomials such that $Q \neq 0$. We make the Euclidean division of
$P$ by $Q$ and we keep 
dividing the intermediate quotients by $Q$ until we arrive at a
quotient of degree $< d(Q)$:

\begin{equation*}
  \begin{cases}
      P = q_{0}Q + r_{0} \\
      q_{0} = q_{1}Q + r_{1} \\
      \vdots \\
      q_{t-1} = q_{t}Q + r_{t}
  \end{cases}.
\end{equation*}

Here $q_{t} \neq 0$ and $d(q_{t}) < d(Q)$. Then we obtain an
expansion of $P$ in terms of $Q$: 
$$ P = q_{t}Q^{t+1} + r_{t}Q^{t} +r_{t-1}Q^{t-1}+\cdots+r_{0}.  $$
All the coefficients $q_{t},r_{t},  r_{t-1},...,r_{0}$ are polynomials
in $Y$ of degrees
$< d(Q)$. This is the unique expansion having this property:

\begin{prop} \label{exp}
   One has a unique 
  \textbf{Q-adic expansion of P}:
  \begin{equation} 
     \label{expansion}
     P = a_{0}Q^{s}+a_{1}Q^{s-1}+\cdots+a_{s}
  \end{equation}
where $a_{0},a_{1},...,a_{s} \in A[Y]$ and $d(a_{i})<d(Q)$ 
for all $i\in \{0,...,s\}$.The $Y$-degrees of the terms $a_{i}Q^{s-i}$
in the right-hand 
side of equation
(\ref{expansion}) are all different 
and $s=[\frac{d(P)}{d(Q)}]$. One has $a_{0}=1$ if and only if
$d(Q) \mid 
d(P)$. In this last situation, supposing that moreover $s$ is
invertible in $A$, one has $a_{1}=0$ if and only if 
$Q = \sqrt[s]{P}$.
\end{prop}

\textbf{Remark:} One should note the analogy with the expansion of
numbers in a basis of numeration. To obtain that notion, one needs
only to take natural numbers in spite of polynomials. Then the $a_{i}$'s 
are the digits of the expansion.

\begin{defin}
 The polynomials $P$ and $Q$ are  given as before, with $d(Q) \mid
 d(P)$. Let us suppose $s=\frac{d(P)}{d(Q)}$ is
 invertible in $A$. The \textbf{Tschirnhausen operator}
 $\tau_{P}$ of  
``completion of the s-power'' is defined by the formula:
\[ \tau_{P}(Q):= Q + \frac{1}{s} a_{1}. \]
\end{defin}

Look again at the example given after Definition \ref{approx}. The
usual expression $P = Y^{n} + \alpha_{1} Y^{n-1} + \cdots + \alpha_{n}$
is the $Y$-adic expansion of $P$ and $\sqrt[n]{P}$ is exactly
$\tau_{P}(Y)$. The following proposition generalizes this observation.

\begin{prop} \label{Iterate}
 Suppose $P \in A[Y]$ is monic and $p \mid d(P)$, with $p$
 invertible in $A$. 
 The approximate roots can be computed by iterating the Tschirnhausen
 operator on arbitrary polynomials of the correct degree:

\[ \sqrt[p]{P} = \underbrace{\tau_{P} \circ \tau_{P} \circ ...\circ
 \tau_{P}}_{d(P)/p} (Q) \]
for all $Q \in A[Y]$ monic of degree $\frac{d(P)}{p}$.
\end{prop}

The steps of the proof of Theorem  \ref{Mainthm} are:

  \textbf{Step 1} Show that there exist polynomials verifying the
  conditions of Theorem  \ref{Mainthm}, point 1).

  \textbf{Step 2} Show that those conditions are preserved by an adequate
  Tschirnhausen operator.

  \textbf{Step 3} Apply Proposition  \ref{Iterate} to show inductively that the
  characteristic roots also satisfy those conditions.

  \textbf{Step 4} Show that the point 2) of  Theorem \ref{Mainthm} is
  true for all polynomials satisfying the conditions of point 1).
\vspace{5mm}

This motivates us to introduce a special name for the polynomials
verifying the conditions of Theorem \ref{Mainthm}, point 1):

\begin{defin} \label{semiroot}
  A polynomial $q_{k} \in \mathbf{C}[[X]][Y]$ is \textbf{a $k$-semiroot} of
  $f$ if it is monic of degree $d(q_{k}) = \frac{N}{E_{k}}$ and 
                $(f, q_{k}) = \overline{B}_{k+1}$.
\end{defin}

The term of ``semiroot'' is taken from \cite{A 89}. 

We show in fact that all the corollaries of the main theorem
(Theorem \ref{Mainthm}), with the exception of the first one, are true
for polynomials that are $k$-semiroots of $f$. That is why we 
begin the proofs of the corollaries \ref{coincid}, \ref{expansions},
\ref{graduation}, \ref{dual} and \ref{curvette} by restating them in this
greater generality. It is only in Corollary \ref{comput}
that the precise construction of approximate roots is useful. In our
context, the
value of the approximate roots lies mainly in the fact that the
definition is global and at the same time gives locally $k$-semiroots (see
section 9). 

We now formulate some propositions that are used
in the proof of Theorem \ref{Mainthm}. The first one is attributed by
some authors to M.Noether. Equivalent statements in terms of characteristic 
exponents can be found in \cite{S 75}, \cite{H 84}, \cite{C 31},
\cite{Z 71}, \cite{M 77}.

\begin{prop} \label{Noether}

 If $\phi \in \mathbf{C}[[X]][Y]$ is monic, irreducible and 
  $K(f,\phi) :=$ $$=\mbox{max} \{ v_{X}(\eta(X)-\zeta(X)), \eta(X)
 \makebox{ and } \zeta(X) 
 \makebox{ are Newton-Puiseux series of } f  \makebox{ and } \phi \} $$
 is \mbox{{\em the coincidence order}} of $f$ and $\phi$, then one
  has the formula: 
 \[ \frac{(f,\phi)}{d(\phi)} = 
 \frac{\overline{B}_{k}}{N_{1}\cdots N_{k-1}} + 
 \frac{N \cdot K(f,\phi)-B_{k}}{N_{1}\cdots N_{k}} \]

 \noindent where $k\in \{0,...,G\}$ is the smallest integer such that 
$K(f,\phi) < \frac{B_{k+1}}{N}$.
\end{prop}

This proposition allows one to translate information about intersection
numbers into information about equalities of truncated Newton-Puiseux series
and conversely. For example, from Definition \ref{semiroot} to Corollary
\ref{coincid}, where in place of $f_{k}$ we consider an arbitrary
semiroot $q_{k}$.

\begin{prop} \label{existence}
  For each $k \in \{0,...,G \}$, the minimal polynomial $\phi_{k}$ of a 
$k$-truncated Newton-Puiseux series $\eta_{k}(X)$ of $f$ is a $k$-semiroot.
\end{prop}

This gives us the Step 1 explained before.

\begin{prop} \label{indep}
  If $\phi \in \mathbf{C}[[X]][Y]$ and $d(\phi)< \frac{N}{E_{k}}$, then 
$(f,\phi) \in \langle \overline{B}_{0},...,\overline{B}_{k} \rangle.$
\end{prop}

In other words, $\frac{N}{E_{k}}$ is the minimal degree for which one
can obtain the value $\overline{B}_{k+1}$ in the semigroup
$\Gamma(C)$.

\begin{prop} \label{stab}
 If $\phi$ is a $k$-semiroot and $\psi$ is a $(k-1)$-semiroot, $k \in
 \{1,...,G \}$, then $\tau_{\phi}(\psi)$ is a $(k-1)$-semiroot of $f$.
\end{prop}

This gives Step 2 in the proof of Theorem \ref{Mainthm}.

\section{The Proofs of the Main Theorem \\and of its Corollaries}

\textbf{Proof of Theorem \ref{Mainthm}} \label{Proofs}

1) The first equality 
$d(f_{k}) = \frac{N}{E_{k}}$ is clear from the definition of
approximate roots. 

The main point is to prove that $(f, f_{k}) = \overline{B}_{k+1}$ 
for all $k \in \{0,...,G \}$, where $\overline{B}_{G+1} = \infty$. We
shall prove it by descending induction, starting from $k=G$. Then
$f_{G}=f$ and so $(f, f_{G}) = \infty = \overline{B}_{G+1}$.

Let us suppose that $(f, f_{k}) = \overline{B}_{k+1}$, with $k \in
\{1,...,G \}$. Then we have by Proposition \ref{tworoots}:
 $$f_{k-1}=\sqrt[E_{k-1}]{f} =
\sqrt[N_{k}E_{k}]{f} = \sqrt[N_{k}]{\sqrt[E_{k}]{f}}$$
and so: $f_{k-1}=\sqrt[N_{k}]{f_{k}}$.

By Proposition \ref{Iterate}, we know that  
 $\sqrt[N_{k}]{f_{k}} = \underbrace{\tau_{f_{k}} \circ \cdots \circ
 \tau_{f_{k}}}_{d(f_{k})/N_{k}} (q_{k-1})$, where $q_{k-1}$ is
 an arbitrary polynomial of degree $\frac{N}{E_{k-1}}$. We shall take
 for $q_{k-1}$ an arbitrary $(k-1)$-semiroot, which exists by
 Proposition \ref{existence}. 
By the induction hypothesis, $ f_{k}$ is a $k$-semiroot. By
 Proposition \ref{stab}, if $\psi$ is a $(k-1)$-semiroot of $f$, then
 $\tau_{f_{k}}(\psi)$ is again a $(k-1)$-semiroot. Starting with
 $\phi_{k-1}$ and applying the operator  $\tau_{f_{k}}$
 consecutively $\frac{d(f_{k})}{N_{k}} = \frac{N}{E_{k-1}}$ times,
 we deduce that $\sqrt[N_{k}]{f_{k}}$ is a $(k-1)$-semiroot of $f$.

The induction step is completed, so we have proved the first part
of the proposition. 

2)We show that this is true generally for an arbitrary
  $k$-semiroot $q_{k}$. First we prove that $q_{k}$ is
  irreducible.

Suppose this is not the case. Then $q_{k}=\prod_{i=1}^{m}r_{i}$, where
$m \geq 2$ and \linebreak $r_{i} \in \mathbf{C}[[X]][Y]$ are monic 
polynomials of
degree at least 1. So, for all $i$, $d(r_{i})< \linebreak <d(q_{k})=
\frac{N}{E_{k}}$. 
By Proposition \ref{indep}, $(f,r_{i}) \in  
\langle \overline{B}_{0},...,\overline{B}_{k} \rangle$ and so 
$(f,q_{k}) = \linebreak =\sum_{i=1}^{m}(f,r_{i}) \in 
\langle \overline{B}_{0},...,\overline{B}_{k} \rangle$, which
contradicts $(f, q_{k}) = \overline{B}_{k+1}$. 
This shows that $q_{k}$ is irreducible. 

We have to prove now the claim concerning its characteristic
exponents. We apply Proposition \ref{Noether}, which expresses 
$\frac{(f, q_{k})}{d(q_{k})}$ in terms of the coincidence
exponent of $f$ and $q_{k}$. 

First, we have directly by the property of being a $k$-semiroot: 
$\frac{(f, q_{k})}{d(q_{k})}=
\frac{\overline{B}_{k+1}}{d(q_{k})}=  \linebreak =
\frac{\overline{B}_{k+1}}{N_{1}\cdots N_{k}}$. So, by Proposition
\ref{Noether}, 
one has $K(f, q_{k}) = B_{k+1}$, which implies that the $k$-truncated
Newton-Puiseux series of $f$ and
$q_{k}$ are equal. This means
that the first $k$ terms of the characteristic sequence of $q_{k}$
are
$l\frac{B_{0}}{E_{k}},l\frac{B_{1}}{E_{k}},...,l\frac{B_{k}}{E_{k}}$,
with $l \in \textbf{N}^*$. So $d(q_{k}) = l\frac{B_{0}}{E_{k}} =
l\frac{N}{E_{k}}$. But we know that $d(q_{k}) = \frac{N}{E_{k}}$,
and this implies that $l=1$, which in turn implies that $q_{k}$ has
no more characteristic exponents. \hfill $\Box$

\medskip

\noindent \textbf{Proof of Corollary \ref{comput}}

The point here is to compute the characteristic approximate roots and
the 
characteristic sequence without previously computing truncated Newton-Puiseux
parameterizations.

The algorithm given in the statement works because 
$$\mbox{gcd}(\overline{B}_{0},...,\overline{B}_{k}) =
\mbox{gcd}(B_{0},...,B_{k}) = E_{k},$$
which is part of Proposition \ref{semigroup}. 

Once the characteristic roots have been computed, by-products of the
algorithm are the sequences $(\overline{B}_{0},...,\overline{B}_{G})$
and $(E_{0},...,E_{G})$. From the point 1 of Proposition \ref{semigroup} one
deduces then the characteristic sequence $(B_{0},...,B_{G})$. \hfill $\Box$

\medskip

 \textbf{Example:}
Take:
$$f(X,Y)=Y^{4}-2X^{3}Y^{2}-4X^{5}Y+X^{6}-X^{7},$$
an example already considered to illustrate Proposition \ref{example}. 
We suppose here we do not know a Newton-Puiseux parameterization for it. We
suppose it is irreducible - indeed it is, and the elaborations of the
algorithm alluded to in the text would show it - so we apply the
algorithm:
$$\begin{array}{l}
  N=B_{0}=\overline{B}_{0}=E_{0}=4 \\
  f_{0}=\sqrt[4]{f}=Y \\
  E_{1}=\mbox{gcd}(4,6)=2 \\
  f_{1}=\sqrt[2]{f}=\tau_{f}\circ \tau_{f}(Y^{2})=Y^{2}-X^{3} \\
  (f,f_{1})=13 \\
  E_{2}=\mbox{gcd}(E_{1},\overline{B}_{2})=\mbox{gcd}(2,13)=1 \\
  G=2 \\
  N_{1}=\frac{E_{1}}{E_{2}}=2 \\
  B_{2}=B_{1}+ \overline{B}_{2} -N_{1}\overline{B}_{1}= 6+13-2\cdot
  6=7.
\end{array}$$
So:
$$(B_{0},B_{1},B_{2})=(4,6,7).$$

\medskip

\noindent \textbf{Proof of Corollary \ref{coincid}}

The more general formulation is:
\medskip
{\it If $q_{k}$ is a $k$-semiroot, $f$ and $q_{k}$ have equal
  $k$-truncated Newton-Puiseux series.} 
\medskip
The proof is contained in that of Theorem 
\ref{Mainthm}, point 2, where it was seen that 
$K(f,q_{k})={B}_{k+1}$. \hfill $\Box$

\medskip

\noindent \textbf{Proof of Corollary \ref{expansions}}

We give first the more general formulation which we prove in the
sequel:
\medskip

{\it
 Let $q_{0},...,q_{G} \in \mathbf{C}[[X]][Y]$ be monic polynomials
 such that for all $i$, $d(q_{i})=\frac{N}{E_{i}}$. Then every 
$\phi \in \mathbf{C}[[X]][Y]$ can be uniquely written in the form:
$$\phi = \sum_{finite}\alpha_{i_{0}...i_{G}}q_{0}^{i_{0}} 
q_{1}^{i_{1}}\cdots q_{G}^{i_{G}} $$
where $i_{G} \in \mathbf{N}$,  $0 \leq i_{k} < N_{k+1}$ for $0 \leq k \leq
G-1$ and the coefficients $\alpha_{i_{0}...i_{G}}$ are elements of the
ring $\mathbf{C}[[X]]$. Moreover:

\hspace{5mm} 1) the $Y$-degrees of the terms appearing in the right-hand
side of the preceding equality are all distinct.

\hspace{5mm} 2) if for every $k \in \{0,...,G\}$, $q_{k}$ is a
$k$-semiroot, then the orders in $T$ of the terms 
$$\alpha_{i_{0}...i_{G}}(T^{N})q_{0}(T^{N},Y(T))^{i_{0}} 
\cdots q_{G-1}(T^{N},Y(T))^{i_{G-1}}$$ 
are all distinct, where $T \rightarrow (T^{N}, Y(T))$ is a Newton-Puiseux
parameterization of $f$.}

\medskip

Take first the $q_{G}$-adic expansion of $\phi$:
$$ \phi = \sum_{0 \leq i_{G} \leq [\frac{d(\phi)}{d(q_{G})}]}\alpha_{i_{G}}
    q_{G}^{i_{G}}. $$

Here $\alpha_{i_{G}} \in \mathbf{C}[[X]][Y]$ and $d(\alpha_{i_{G}})<
 d(q_{G}) = \frac{N}{E_{G}}=N.$

Take now the $ q_{G-1}$-adic expansion of every coefficient
$\alpha_{i_{G}}$:
$$ \alpha_{i_{G}} = \sum \alpha_{i_{G-1}i_{G}}q_{G-1}^{i_{G-1}}. $$
The coefficients $\alpha_{i_{G-1}i_{G}} \in \mathbf{C}[[X]][Y]$ have
degrees $d(\alpha_{i_{G-1}i_{G}}) < d(q_{G-1})$ and the sum is over 
$i_{G-1} < N_{G}$.

Proceeding in this manner we get an expansion with the required properties.
Before proving the unicity, we prove point 1), namely the inequality of
the degrees. 

Suppose there exist $(i_{0},...,i_{G}) \neq (j_{0},...,j_{G})$ and 
$$d(\alpha_{i_{0}...i_{G}}q_{0}^{i_{0}} 
q_{1}^{i_{1}}\cdots q_{G}^{i_{G}}) = 
d(\alpha_{j_{0}...j_{G}}q_{0}^{j_{0}} 
q_{1}^{j_{1}}\cdots q_{G}^{j_{G}}) \neq \infty.$$

This means:
$$\sum_{k=0}^{G}i_{k}\cdot \frac{N}{E_{k}} = 
\sum_{k=0}^{G}j_{k}\cdot \frac{N}{E_{k}}. $$

Let us define $p \in \{0,...,G\}$ such that $i_{k}=j_{k}$ for $k\geq
p+1$ and $i_{p} < j_{p}$. If such a $p$ does not exist, simply interchange 
$(i_{0},...,i_{G})$ and $(j_{0},...,j_{G})$, then apply the preceding
definition. We obtain: 
$$ \sum_{k=0}^{p-1}(i_{k}-j_{k})\frac{N}{E_{k}} = 
 (j_{p}-i_{p})\frac{N}{E_{p}}. $$

But $j_{p}-i_{p} \geq 1$ and $\mid i_{k}-j_{k} \mid \leq N_{k+1}-1$,
so:
\begin{equation*}
  \begin{array}{ll}
      \frac{N}{E_{p}} & \leq
      \sum_{k=0}^{p-1}(N_{k+1}-1)\frac{N}{E_{k}} =
      \sum_{k=0}^{p-1}(\frac{E_{k}}{E_{k+1}}-1)\frac{N}{E_{k}} = \\
                      & =
      \sum_{k=0}^{p-1}(\frac{N}{E_{k+1}}-\frac{N}{E_{k}}) =
      \frac{N}{E_{p}}-1
  \end{array}
\end{equation*}
\noindent which is a contradiction.

Now, this property of the degrees shows that $0 \in
\mathbf{C}[[X]][Y]$ has only the trivial expansion, and this in turn
shows the unicity.

Let us move to point 2). From now on, $q_{k}$ is a $k$-semiroot. 
By the properties of intersection numbers recalled in
section 4, $v_{T}(q_{k}(T^{N},Y(T))) =
(f,q_{k})=\overline{B}_{k+1}.$ So:
$$ v_{T}(\alpha_{i_{0}...i_{G}}(T^{N})q_{0}(T^{N},Y(T))^{i_{0}} 
\cdots q_{G-1}(T^{N},Y(T))^{i_{G-1}})=
\sum_{k=-1}^{G-1}i_{k}\overline{B}_{k+1}.$$
Here $i_{-1}=v_{X}(\alpha_{i_{0}...i_{G}}(X)) \in \mathbf{N}$.

Let us suppose we have $(i_{0},...,i_{G-1}) \neq (j_{0},...,j_{G-1})$
such that: $\sum_{k=-1}^{G-1}i_{k}\overline{B}_{k+1}= \linebreak =
\sum_{k=-1}^{G-1}j_{k}\overline{B}_{k+1}$. As before, we take
$p \in \{0,...,G-1\}$ with $i_{k}=j_{k}$ for $k \geq p+1$ and $i_{p} <
j_{p}$. So: $(j_{p}-i_{p})\overline{B}_{p+1} = 
\sum_{k=-1}^{p-1}(i_{k}-j_{k})\overline{B}_{k+1}$ which gives:
$E_{p} \mid (j_{p}-i_{p})\overline{B}_{p+1}.$ But $E_{p+1}= 
\mbox{gcd}(E_{p}, \overline{B}_{p+1})$, by Proposition
\ref{semigroup}, and we get:  
$N_{p+1}=\frac{E_{p}}{E_{p+1}} \mid (j_{p}-i_{p})$. As $0 <
j_{p}-i_{p} <N_{p+1}$, we get a contradiction.

With this, point 2) is proved.\hfill $\Box$

\medskip

\noindent \textbf{Proof of Corollary \ref{graduation}}

We prove the following fact:
\medskip

{\it If $q_{0},...,q_{G}$ are semiroots of $f$, the images of 
$(X,q_{0},...,q_{G-1})$ in the graded ring
$gr_{v_{T}}(\mathcal{O}_{C})$ generate it. If the coordinates are
generic, they form a minimal system of generators. }

\medskip

We take the notations explained in section 5, with
$A=\mathcal{O}_{C}$. For every $p \in \linebreak \in\Gamma(C), 
\mbox{ dim}_{\mathbf{C}}(I_{p}/I_{p}^{+})=1$. The vector space
$I_{p}/I_{p}^{+}$ is generated by an arbitrary element $\phi \in
\mathcal{O}_{C}$ such that $v_{T}(\phi)=p$. We obtain: 
$$ gr_{v_{T}}(\mathcal{O}_{C}) \simeq \bigoplus_{\{p \in \Gamma(C)\}}
\mathbf{C}T^{p}.$$
 We have: $v_{T}(X)=N$ and $v_{T}(q_{k})=\overline{B}_{k+1}$, for $k
 \in \{0,...,G-1\}$. To show that the images of $X,q_{0},...,q_{G-1}$
 generate $gr_{v_{T}}(\mathcal{O}_{C})$ is equivalent with the fact
 that every $\omega \in gr_{v_{T}}(\mathcal{O}_{C})$ can be expressed
 as a polynomial $P_{\omega}(T^{N},T^{\overline{B}_{1}},...,
 T^{\overline{B}_{G}})$. This comes in turn from Proposition
 \ref{semigroup}. Indeed, it is 
 shown that $\langle
 \overline{B}_{0},\overline{B}_{1},...,\overline{B}_{G} \rangle =
\linebreak =
 \Gamma(C)$ and so every $p \in \Gamma(C)$ can be written as $p = 
 \sum_{k=-1}^{G-1}i_{k}\overline{B}_{k+1}$, which implies
 $T^{p}=(T^{N})^{i_{-1}} (T^{\overline{B}_{1}})^{i_{0}}\cdots 
(T^{\overline{B}_{G}})^{i_{G-1}}$. An arbitrary 
$\omega \in gr_{v_{T}}(\mathcal{O}_{C})$ is then a linear combination
of such terms.

Another proof can use Corollary \ref{expansions}.

In case the coordinates are generic, $\overline{B}_{0}=m(C)$, the
multiplicity of $C$ at the origin, and this is the smallest non-zero
value in $\Gamma(C)$. Then
$(\overline{b}_{0},...,\overline{b}_{g})$ is a minimal system
of generators of $\Gamma(C)$. Indeed, what prevented 
$(\overline{B}_{0},...,\overline{B}_{G})$ from being minimal was the
possibly non minimal value of $\overline{B}_{0}$ in $\Gamma(C)-\{0\}$
(see Proposition \ref{semigroup}).

Now, the minimality for the algebra $gr_{v_{T}}(\mathcal{O}_{C})$
comes from the minimality for the semigroup $\Gamma(C)$.\hfill $\Box$

\textbf{Remark:} An equivalent statement (using the notion of maximal
contact explained after Corollary \ref{dual}, rather than the notion
of semiroot), was proved by M. Lejeune-Jalabert. See the paragraph
1.2.3 in the Appendix of \cite{Z 73}.

\medskip

\noindent \textbf{Proof of Corollary \ref{dual}}

Instead of the characteristic roots we consider arbitrary semiroots
$q_{k}$ and we show that the Corollary is also true in this greater
generality. We sketch three proofs of the Corollary. The first one
uses adequate coordinate systems to follow the strict transforms of
$C_{f}$ and $C_{q_{k}}$ during the process of blowing-ups. The second
and third one are more intrinsic.

1) 
Let us consider generic coordinates $(X,Y)$ and Newton-Puiseux series
$\eta(X),$ \linebreak  $\eta_{k}(X)$ for $f$, respectively $q_{k}$. We have:
$$ \eta(X) = \sum_{j \geq n}a_{j}X^{\frac{j}{n}}.$$

If $\gamma_{1}: S_{1} \rightarrow \mathbf{C}^{2}$ is the blow-up
of $0 \in \mathbf{C}^{2}$, the strict transform $C_{f}^{1}$ of $C_{f}$
in $S_{1}$ 
passes through the origin of a chart of coordinates $(X_{1},Y_{1})$ such
that:
$$\begin{cases}
    X = X_{1} \\
    Y = X_{1}(a_{n}+ Y_{1})
  \end{cases}.$$

The strict transform $C_{f}^{1}$ of $C_{f}$ has in the coordinates 
$(X_{1},Y_{1})$ a Newton-Puiseux series of the form:
$$ \eta^{(1)}(X_{1}) = \sum_{j \geq n+1}a_{j}X_{1}^{\frac{j}{n}-1}.$$
     
The coordinates $(X_{1},Y_{1})$ are generic for it if and only if 
$[\frac{b_{1}}{n}] \geq 2$. If this is the case, one describes the
restriction of the next blowing-up to the chart containing the strict
transform of $C_{f}$ by the change of variables:
$$\begin{cases}
    X_{1} = X_{2} \\
    Y_{1} = X_{2}(a_{2n}+ Y_{2})
  \end{cases}.$$

One continues like this $s_{1}:=[\frac{b_{1}}{n}]$ times till one
arrives in the chart $(X_{s_{1}},Y_{s_{1}})$ at a strict transform 
$C_{f}^{s_{1}}$ with Newton-Puiseux series:
$$ \eta^{(s_{1})}(X_{s_{1}}) = 
\sum_{j \geq b_{1}}a_{j}X_{s_{1}}^{\frac{j}{n}-s_{1}}.$$

Now for the first time the coordinates are not generic with
respect to the 
series. Let us look also at the strict transform
$C_{q_{0}}^{s_{1}}$ of $C_{q_{0}}$. By Corollary \ref{coincid}, the
branch $C_{q_{0}}$ has  a Newton-Puiseux series $\eta_{0}(X)$ such that:
$$ \eta_{0}(X) = \sum_{j \geq 1}a_{j}'X^{\frac{j}{n}} $$
with $a_{j}' = a_{j}$ for $j < b_{1}$ and $n \mid j$ for all $j \in
\mathbf{N}^*$. 

The strict transform $C_{q_{0}}^{s_{1}}$ then has a Newton-Puiseux series of
the form:
$$  \eta_{0}^{s_{1}}(X_{s_{1}}) = 
\sum_{j \geq b_{1}}a_{j}'X_{s_{1}}^{\frac{j}{n}-s_{1}}. $$

The series in the right-hand side has integral exponents, which shows that 
$C_{q_{0}}^{s_{1}}$ is smooth - which was evident, as $C_{q_{0}}$ was
already smooth. But, more important, $C_{q_{0}}^{s_{1}}$ is not tangent to
$X_{s_{1}}=0$. This shows that it is transverse to $C_{f}^{s_{1}}$
and to the only component of the exceptional divisor passing through 
$(X_{s_{1}},Y_{s_{1}})=(0,0)$, which is defined by the equation:
$X_{s_{1}}=0$.

The next blowing-up separates the strict transforms of $C_{f}$ and
$C_{q_{0}}$. The curve $C_{q_{0}}^{(s_{1}+1)}$ passes through a smooth
point of the newly created component of the exceptional divisor. 

\begin{figure}
\begin{center}
\includegraphics[width=3.5in]{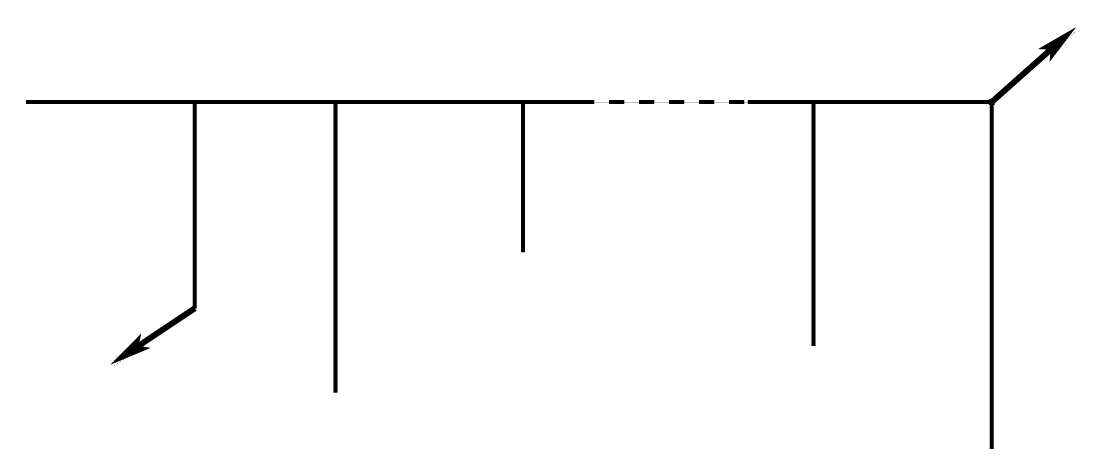}
\end{center}
\caption{The Dual Graph of the Product}
\label{proddualgr}

\end{figure}

This shows that the dual graph of the total transform of $f \cdot
f_{0}$ is as drawn in Figure \ref{proddualgr}.

To continue, one needs to change coordinates after $s_{1}$
blowing-ups. Instead of considering the ordered coordinates 
$(X_{s_{1}},Y_{s_{1}})$, we look at $(Y_{s_{1}},X_{s_{1}})$. We now
use 
the inversion formulae explained in Proposition \ref{Inversion}. They
allow to express the characteristic exponents of $C_{f}^{s_{1}}$ with
respect 
to $(Y_{s_{1}},X_{s_{1}})$ in terms of those with respect to 
$(X_{s_{1}},Y_{s_{1}})$. Moreover, it follows from the property of
truncations stated 
in Proposition \ref{Inversion} that, if one inverts
simultaneously the strict transforms of $C_{f},
C_{q_{0}},...,C_{q_{g-1}}$, they keep having coinciding Newton-Puiseux series
up to controlled orders. Repeating this process, one shows that after
a number of inversions equal to the number of terms in the continuous
fraction expansion of $\frac{b_{1}}{n}$, the initial situation
is repeated, but with a curve having genus $(g-1)$. The strict
transform of the semiroot $q_{k}$, for $k \in
\{1,...,g\}$ will be a $(k-1)$-semiroot for the strict transform of
$f$, in the natural
coordinates resulting from the process of blowing-ups. So, one can
iterate the analysis made for $q_{0}$ and get the corollary. 
\hfill $\Box$

\medskip

2)
Given two branches at the origin, from the knowledge of their characteristic 
exponents and of their coincidence exponent (see
Proposition \ref{Noether}), one can construct the dual graph of resolution of
their product. This is explained in \cite{LMW 89} and proved in
detail, as well as in the case of an arbitrary number of branches, in
\cite{GB 96}. In our case this shows that the minimal embedded
resolution of $f$, where $q_{k}$ is an arbitrary semiroot
for generic coordinates, is also an embedded resolution of $f \cdot
q_{k}$. Moreover, the extended dual graph is obtained from the dual
graph $D(\pi_{m},f)$ attaching an arrow-head vertex at the end of the
$k$-th vertical segment. (see the explanations given after Corollary 
\ref{dual}). We get from it the corollary.
\hfill $\Box$

\medskip

3)
If $l \in \mathbf{C}[[X,Y]]$ is of multiplicity $1$, let
$\pi_{m}^*(l)$ be its total transform divisor on $\Sigma_{m}$. If $L$
is a component of $\mathcal{E}$, the exceptional divisor of $\pi_{m}$, let
$\mu(L)$ be its multiplicity in $\pi_{m}^*(l)$. Let also $\phi(L)$ be
its multiplicity in $\pi_{m}^*(f)$. These multiplicities can be
computed inductively, following the order of creation of the
components in the process of blowing-ups. In particular, if $L_{k}$ is
the component represented at the end of the $k$-th vertical segment of
the dual graph, $\mu(L_{k})= \frac{n}{e_{k-1}}$ for $k \in \{1,...,g\}$,
and $\phi(L_{k})= \overline{b}_{k}$ (folklore). 

Let us consider now the branch $C_{q_{k-1}}$. We know that $(f,
q_{k-1}) = \overline{b}_{k}$ and $m(q_{k-1})=\frac{n}{e_{k-1}}$,
where $n = m(f)$. Then $\frac{(f,q_{k-1})}{ m(q_{k-1})} =
\frac{e_{k-1}\overline{b}_{k}}{n}$ and the lemma on the growth of
coefficients of insertion in \cite{LMW 89} shows that the strict
transform of $C_{q_{k-1}}$ necessarily meets a component of the $k$-th
vertical segment of $D(\pi_{m})$. If $C_{q_{k-1}}'$ is the strict
transform of $q_{k-1}$ by $\pi_{m}$, we have:  
$$ (f,q_{k-1}) = ( \pi_{m}^*(f), C_{q_{k-1}}') = 
\sum_{L}\phi(L)(L,C_{q_{k-1}}')$$
the sum being taken over all the components of $\mathcal{E}$ which
meet 
$C_{q_{k-1}}'$. Now it can be easily seen that $\phi$ strictly grows
on a vertical segment of $D(\pi_{m})$, from the end to the point of
contact with the horizontal segment. This comes from the fact that
those components of the exceptional divisor are created in this order
- but not necessarily consecutively. As $\phi(L_{k})=
\overline{b}_{k}$ and $(f, 
q_{k-1}) = \overline{b}_{k}$, we see that: 
$$\overline{b}_{k} = 
\sum_{L}\phi(L)(L,C_{q_{k-1}}') \geq \sum_{L}\phi(L)m(C_{q_{k-1}}') 
\geq \phi(L_{k})\cdot 1 = \overline{b}_{k}.$$
This means that the inequalities are in fact equalities and shows that 
$C_{q_{k-1}}'$ is smooth, meets $L_{k}$ transversely and meets no other
component of $\mathcal{E}$. 
\hfill $\Box$

\medskip

\noindent \textbf{Proof of Corollary \ref{curvette}}

We prove:
\medskip

\textit{A semiroot $q_{k}$ of $f$ in arbitrary coordinates is a
  curvette with respect to one of the components $L_{k},L_{k+1}$.}
\medskip

We analyze successively the three cases introduced in
Proposition \ref{Inversion}, using also some results of its proof.
 \medskip

1) $\underline{B_{0}=b_{0}}.$

This is the case of generic coordinates. The affirmation is the same
as Corollary \ref{dual}. We get that $q_{k}$ is a curvette with
respect to $L_{k+1}$, for all $k \in \{0,...,g\}$.
\medskip

2) $\underline{B_{0}=lb_{0}},$ with $2 \leq l \leq
   [\frac{b_{1}}{b_{0}}] .$

Then, by Proposition \ref{Inversion}, $G=g+1$.

The curve $q_{0}$ is smooth and so $m(q_{0})=1$.

Moreover, by the definition of semiroots, $(f,q_{0})=\overline{B}_{1}=
\overline{b}_{0}=m(f).$

This shows that $q_{0}$ is smooth and transversal to $f$ and so it is
a curvette with respect to $L_{0}$. 

If $k \in\{1,...,G\}$, where by
Proposition \ref{Inversion}, $G=g+1$, we have:
$$ m(q_{k})=b_{0}(q_{k})=B_{1}(q_{k})=\frac{B_{1}}{E_{k}}=
\frac{b_{0}}{e_{k-1}},$$
$$(f,q_{k})=\overline{B}_{k+1}=\overline{b}_{k}.$$ 
We have noted by $b_{0}(q_{k})$ the corresponding characteristic exponent
of $q_{k}$.

So, for $k \in \{1,...,G-1\}$, the curve $q_{k}$ is a $(k-1)$-semiroot
with respect to generic
coordinates and by Corollary \ref{dual}, it is a curvette with respect to
$L_{k}$.
\medskip

3) $\underline{B_{0}=b_{1}}$

By Proposition \ref{Inversion}, we have $G=g$.

Again $q_{0}$ is smooth. Using Proposition \ref{Inversion} we obtain:
$$ (f,q_{0})=\overline{B}_{1}=\overline{b}_{0}=m(f).$$
As in the preceding case, $q_{0}$ is a curvette with respect to
$L_{0}$.

If $k\in\{1,...,G\}$, where $G=g$, we have:
$$m(q_{k})=b_{0}(q_{k})=B_{1}(q_{k})=\frac{B_{1}}{E_{k}}=
\frac{b_{0}}{e_{k}},$$
$$(f,q_{k})=\overline{B}_{k+1}=\overline{b}_{k+1}.$$

So $q_{k}$ is a $k$-semiroot with respect to generic coordinates, and
by Corollary \ref{dual}, it is a curvette with respect to $L_{k+1}$.

\begin{figure}
\begin{center}
\includegraphics[width=3.5in]{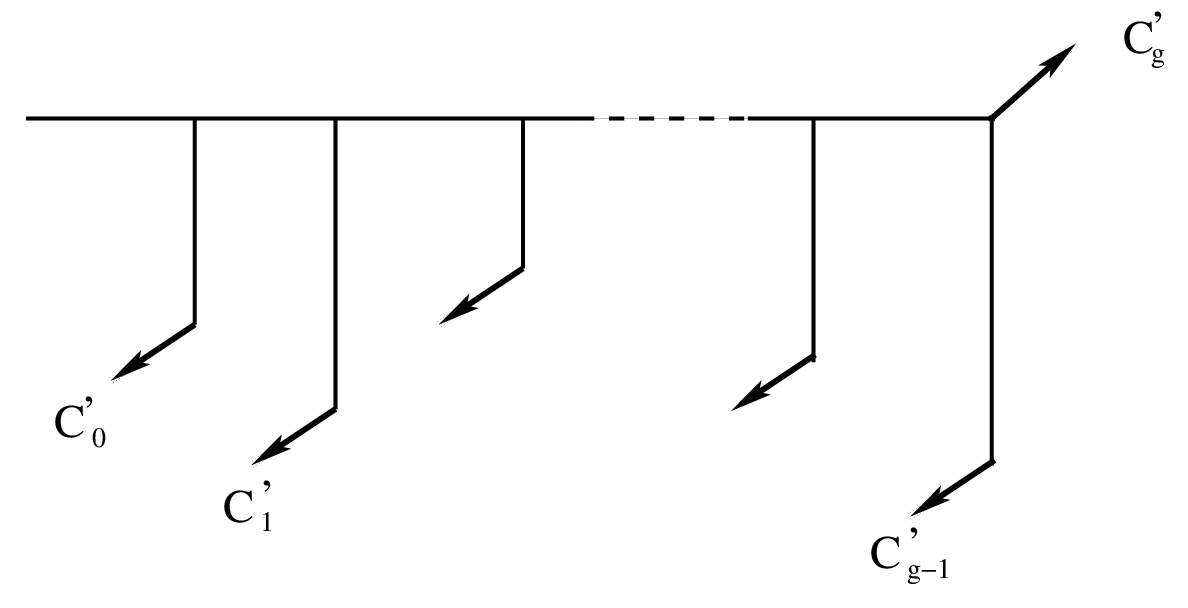}
\end{center}
\caption{$B_{0}= b_{0}$}
\label{case1}
\end{figure}

\begin{figure}
\begin{center}
\includegraphics[width=3.5in]{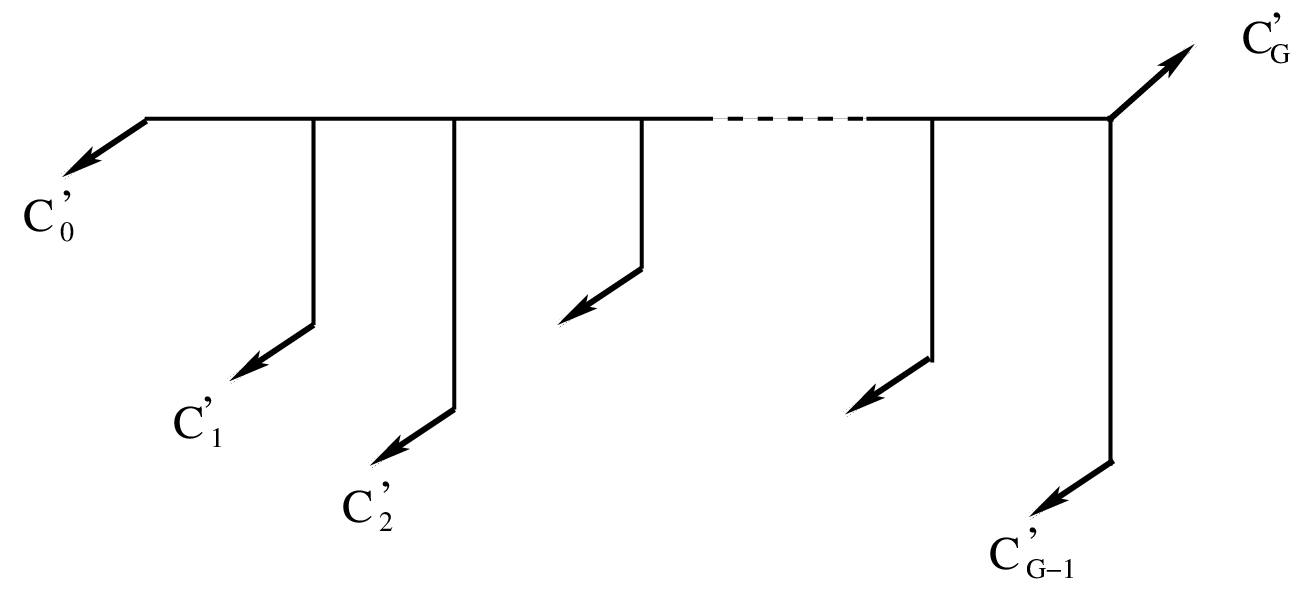}
\end{center}
\caption{$B_{0}= l b_{0}, l \geq 2$}
\label{case2}
\end{figure}

\begin{figure}
\begin{center}
\includegraphics[width=3.5in]{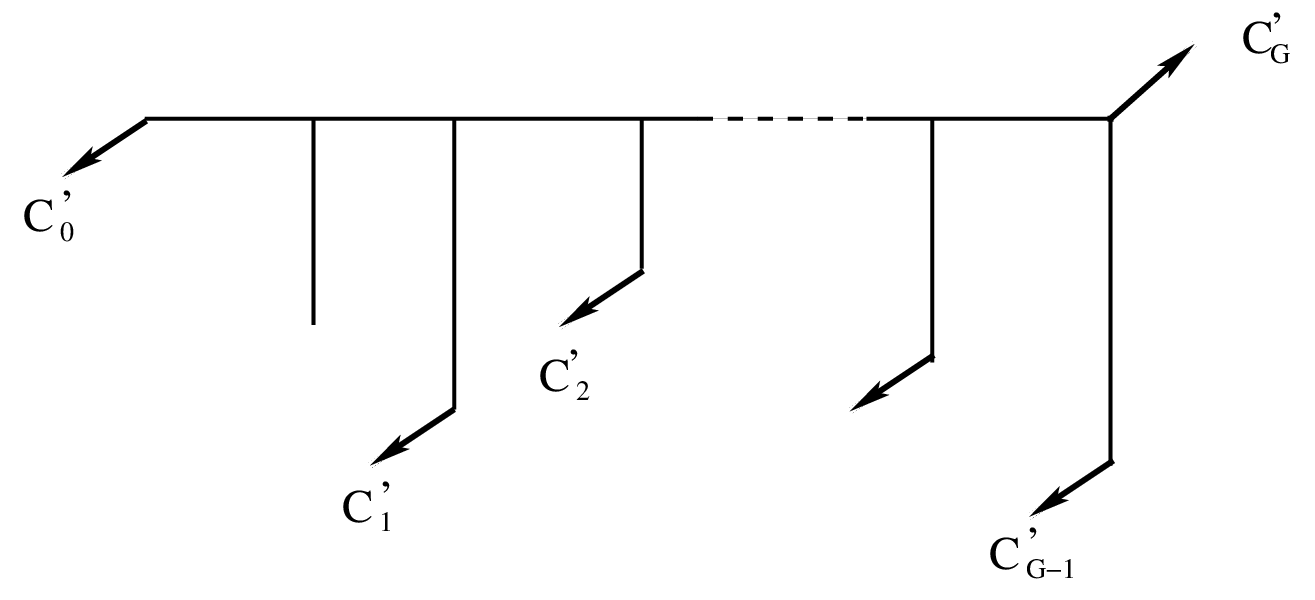}
\end{center}
\caption{$B_{0}= b_{1}$}
\label{case3}
\end{figure}

Let us summarize this study by drawing for each of the three cases the
dual graph of the total transform of the product $q_{0}\cdots
q_{G}$. As in the statement of Corollary \ref{dual}, we denote by
$C_{k}'$ the strict transform of $q_{k}$. We obtain the situations
indicated in  Figures \ref{case1}, \ref{case2}, \ref{case3}.

\hfill $\Box$

\section{The Proofs of the Propositions}

\noindent\textbf{Proof of Proposition \ref{differences}}

The series $\zeta(X)$ can be obtained from $\eta(X)$ by replacing
$X^{\frac{1}{N}}$ by $\omega X^{\frac{1}{N}}$, where $\omega \in
\mu_{N}$, the group
of $N$-th roots of unity. One has the inclusions of cyclic groups:
$\mu_{\frac{N}{E_{0}}} \subset \mu_{\frac{N}{E_{1}}} \subset \cdots \subset 
\mu_{\frac{N}{E_{G}}}=\mu_{N}.$ Let $k \in\{0,...,G\}$ be such that
$\omega \in \mu_{\frac{N}{E_{k}}}-\mu_{\frac{N}{E_{k-1}}}$. Then: 
$$v_{X}(\eta(X)-\zeta(X))= \begin{cases}
                             \infty, \mbox{ if } k=0 \\
                             \frac{E_{k}}{N}, \mbox{ if } k \in
                             \{1,...,G\}
                           \end{cases}.$$
\hfill $\Box$

\noindent\textbf{Proof of Proposition \ref{example}}

We start from $f = Y^{N} + \alpha_{1}(X)Y^{N-1} +
\alpha_{2}(X)Y^{N-2}+\cdots +
\alpha_{N}(X)$. Let us consider the approximate root 
$f_{k} = \sqrt[E_{k}]{f}$.

As is seen from equation (\ref{recursion}) in the proof of Proposition
\ref{defrt}, its coefficients depend only on  
$\alpha_{1}(X),...,\alpha_{\frac{N}{E_{k}}}(X)$. 

Corollary \ref{coincid} shows that $f$ and $f_{k}$ have equal $k$-truncated
Newton-Puiseux series. Combining these facts we see that the $k$-truncated
Newton-Puiseux series of $f$ depend only on $\alpha_{1}(X),...,
\alpha_{\frac{N}{E_{k}}}(X)$. \hfill $\Box$

\medskip

\noindent \textbf{Proof of Proposition \ref{defrt}}

Let us put 
$$ Q =  Y^{\frac{n}{p}}+ a_{1}Y^{\frac{n}{p}-1} +\cdots+ 
a_{\frac{n}{p}}. $$

The inequality  $d(P -Q^{p}) < d(P)-\frac{d(P)}{p}$ means that the
coefficients of $Y^{n},$ \linebreak $Y^{n-1},...,Y^{n-\frac{n}{p}}$ in the
polynomial $P - Q^{p}$ are equal to 0. This gives the system of
equalities:

\begin{equation} 
\label{recursion}
 \begin{cases}
    \alpha_{1} = p a_{1}\\
    \alpha_{2} = p a_{2} + \left( \begin{array}{l}
                                       p \\ 2
                                   \end{array} \right) a_{1}^{2}\\
    \vdots \\
    \alpha_{k} = p a_{k} + \sum_{i_{1}+2i_{2}+\cdots+(k-1)i_{k-1}=k}
c_{i_{1}...i_{k-1}}a_{1}^{i_{1}}\cdots a_{k-1}^{i_{k-1}}, \; 1 \leq k \leq 
\frac{n}{p}
 \end{cases}   
\end{equation}
\vspace{5mm}

Here the coefficients $c_{i_{1}...i_{k-1}}$ are integers, easily
expressible in terms of binomial coefficients:

\[ c_{i_{1}...i_{k-1}} = \left( \begin{array}{c}
                                       p \\ i_{1}+\cdots +i_{k-1}
                                   \end{array} \right)
\frac{(i_{1}+\cdots +i_{k-1})!}{i_{1}!\cdots i_{k-1}!}. \]

 We see that
from the relations (\ref{recursion}) one can compute successively 
$a_{1},a_{2},...,a_{\frac{n}{p}}$. One has only to divide at each step
by $p$. That is the reason why in the
definition of the approximate root we asked $p$ to be invertible.

So $a_{1},a_{2},...,a_{\frac{n}{p}}$ exist and are uniquely
determined. Moreover, they depend only on 
$\alpha_{1},..., \alpha_{\frac{n}{p}}$. \hfill $\Box$

\medskip

\noindent \textbf{Proof of Proposition \ref{tworoots}}

Let us note $Q :=\sqrt[p]{P}$ and $R:=\sqrt[q]{Q}$.

We want to show that $R=\sqrt[pq]{P}$, i.e. that 
$d(P -R^{pq}) < d(P)-\frac{d(P)}{pq}$.

If $S:=Q-R^{q}$, we know that $d(S)<d(Q)-\frac{d(Q)}{q}= 
\frac{d(P)}{p}-\frac{d(P)}{pq}.$ Then: $P-Q^{p}= P-(R^{q}+S)^{p} = 
(P-R^{pq})-\sum_{k=1}^{p}\left( \begin{array}{l}
                                       p \\ k
                                   \end{array}
                                 \right)S^{k}R^{q(p-k)}$, 
and so:
 
$$P-R^{pq}=(P-Q^{p})+\sum_{k=1}^{p}\left( \begin{array}{l}
                                           p \\ k
                                          \end{array}
                                               \right)S^{k}R^{q(p-k)}$$
which implies:

$$ d(P-R^{pq})\leq \mbox{max}(\{d(P-Q^{p})\}\cup \{d(S^{k}R^{q(p-k)}), 1 \leq
k \leq p \}). $$

We know that $d(P-Q^{p}) <d(P)-\frac{d(P)}{p}$, and for $1 \leq k \leq
p$ we have:

\begin{eqnarray*}
   d(S^{k}R^{q(p-k)}) & = & kd(S) + q(p-k)d(R) < \\
                      & < & k(\frac{d(P)}{p}-\frac{d(P)}{pq})+ 
                                    q(p-k)\frac{d(P)}{pq} = \\
                      & = & k\frac{d(P)}{p}-k\frac{d(P)}{pq}+d(P)-  
                                    k\frac{d(P)}{p} = \\
                      & = & d(P)- k\frac{d(P)}{pq} \leq
                                     d(P)-\frac{d(P)}{pq}.
\end{eqnarray*}

So finally: 
$$ d(P-R^{pq})< \mbox{max}\{d(P)-\frac{d(P)}{p}, d(P)-\frac{d(P)}{pq}\} = 
                d(P)-\frac{d(P)}{pq} $$
which shows that $R = \sqrt[pq]{P}.$ \hfill $\Box$

\medskip

\noindent \textbf{Proof of Proposition \ref{second}}

If $P(Y) = Y^{n} + \alpha_{1}Y^{n-1} + \alpha_{2}Y^{n-2}+\cdots +
\alpha_{n}$, then :
$$ P_{1}(Z) =
Z^{-n}(1+\alpha_{1}Z+\cdots +\alpha_{n}Z^{n})$$ 
and so:
$$ P_{1}^{\frac{1}{p}}(Z) = Z^{-\frac{n}{p}}(1 + \sum_{k \geq 1}
c_{k}Z^{k}) = M(P_{1}^{\frac{1}{p}}) + H(P_{1}^{\frac{1}{p}}) $$
where:
$$ M(P_{1}^{\frac{1}{p}}): = Z^{-\frac{n}{p}} + c_{1}
Z^{1-\frac{n}{p}}+\cdots +c_{\frac{n}{p}}, $$
$$ H(P_{1}^{\frac{1}{p}}): = \sum_{k \geq 1}c_{k+\frac{n}{p}}Z^{k}. $$

Here the coefficients $c_{k}$ are elements of $A$, uniquely determined
polynomially by the coefficients of $P$.
We get:
$$ Q(Y) =  Y^{\frac{n}{p}}+ c_{1}Y^{\frac{n}{p}-1} +\cdots + 
c_{\frac{n}{p}}. $$

Let us consider:
$$ R(Y) := P(Y) -Q(Y)^{p}. $$
We want to show that $d(R) < n - \frac{n}{p}$, which is equivalent to 
$v_{Z}(R(Z^{-1})) \geq -n + \frac{n}{p} +1$, $v_{Z}$ designating the
order of a series in $A((Z))$. But:
\begin{equation*}
 \begin{array}{ll}
  R(Z^{-1}) & = P(Z^{-1})- Q(Z^{-1})^{p} = \\
            & = P_{1}(Z) - (M(P_{1}^{\frac{1}{p}}))^{p} = \\
            & = P_{1} - (P_{1}^{\frac{1}{p}} -
            H(P_{1}^{\frac{1}{p}}))^{p} = \\
            &  = \sum_{k=1}^{p}(-1)^{k+1}\left( \begin{array}{l}
                                       p \\ k
                                   \end{array} \right) 
                 P_{1}^{\frac{p-k}{p}}S^{k}.
 \end{array}
\end{equation*}

\noindent where we have noted $S := H(P_{1}^{\frac{1}{p}})$.
We obtain:
$$v_{Z}(R(Z^{-1})) \geq \mbox{min}_{1\leq k \leq p}
    \{v_{Z}(P_{1}^{\frac{p-k}{p}}S^{k})\}  =
   \mbox{min}_{1\leq k \leq p}\{ -n\cdot \frac{p-k}{p} + k\cdot 1 \} = 
   -n+\frac{n}{p} +1$$ 
which is the inequality we wanted to prove. 

So $d(P(Y) -Q(Y)^{p}) < n - \frac{n}{p}$, and this shows that $Q =
\sqrt[p]{P}$. \hfill $\Box$

\medskip

\noindent \textbf{Proof of Proposition \ref{semigroup}}

The degree $d(f)$ can be obtained as an intersection number: $N=d(f)= 
v_{T}(f(0,T))=(f,X)$. So $N \in \Gamma(C)$.

We now {\it define} $\overline{b}_{k}$ by the relation given in 
point 2) of the proposition. We  prove that the numbers defined
in this way are indeed elements of the semigroup $ \Gamma(C)$
and verify the minimality property used to define them in the text.

The important fact is that Proposition \ref{Noether} is proved only using the
formulas of the $\overline{B}_{k}$'s in terms of the $B_{k}$'s. 
That is the reason why we can
apply it in what follows. 

Consider the minimal polynomials $\phi_{k}$ of the $k$-truncated
Newton-Puiseux series $Y_{k}(X)$, for $k \in \{0,...,G-1\}$ (see
Proposition \ref{existence} and its proof). Then $d(\phi_{k})=
\linebreak =
\frac{N}{E_{k}} = N_{1}\cdots N_{k}$ and $K(f, \phi_{k}) =
\frac{B_{k+1}}{N}$, so Proposition \ref{Noether} gives: 
$\frac{(f, \phi_{k})}{d(\phi_{k})} = 
\frac{\overline{B}_{k+1}}{N_{1}\cdots N_{k}}$. We get:
$$ (f, \phi_{k})= \overline{B}_{k}.$$

This shows that $\overline{B}_{k} \in \Gamma(C)$.

Consider now an arbitrary element $g \in \mathbf{C}[[X]][Y]$ and
expand it in terms of $(\phi_{0},...,\phi_{G})$ as explained in the
proof of 
Corollary \ref{expansions}. Indeed, $(\phi_{0},...,\phi_{G})$ are
semiroots of $f$ and we show in the proof that the corollary is true
in this greater generality. Notice that the content of this corollary
is true in our case, because we use only point 2) of
Proposition  \ref{semigroup} which, as well as point 3), results from
point 1).  
 
From Corollary \ref{expansions} we get:
   \begin{equation*}
     \begin{array}{ll}
       (f,g) & = \mbox{min}_{(i_{0},...,i_{G})}
                   \{
                  v_{T}( \alpha_{i_{0}...i_{G}}(T^{N})
                          \phi_{0}(T^{N},Y(T))^{i_{0}}
                                \cdots \phi_{G-1}(T^{N},Y(T))^{i_{G-1}})
                   \} = \\
             & = \mbox{min}_{(i_{0},...,i_{G})}
                    \{i_{-1}N+i_{0}\overline{B}_{1}+\cdots +
                       i_{G-1}\overline{B}_{G} \}
     \end{array}
   \end{equation*}
where $ i_{-1} = v_{X}( \alpha_{i_{0}...i_{G}})$.

This shows that $ \Gamma(C) = 
\langle \overline{B}_{0},...,\overline{B}_{G} \rangle$.     

Now, for every $k \in \{1,...,G \}$, we have $ \overline{B}_{k} \notin 
\langle \overline{B}_{0},...,\overline{B}_{k-1} \rangle$, because
$E_{k-1}$ does not divide $\overline{B}_{k}$.

Suppose $l \in \Gamma(C)$ and $l \notin 
\langle \overline{B}_{0},...,\overline{B}_{k-1} \rangle$. We already
 know that $l \in 
 \langle \overline{B}_{0},...,\overline{B}_{G} \rangle$, so 
$l = i_{-1}n+i_{0}\overline{B}_{1}+\cdots +
                       i_{G-1}\overline{B}_{G}$ with $i_{j} \in
                       \mathbf{N}$ for $j \in \{-1,...,G-1 \}.$
As $l \notin \langle \overline{B}_{0},...,\overline{B}_{k-1} \rangle$,
we deduce that for some $j \geq k$, we have  $i_{j} > 0,$ which 
implies $l \geq 
\overline{B}_{j} \geq \overline{B}_{k}$. This proves the equality 
we were seeking:
\[ \overline{B}_{k} = \mbox{min}\{j \in \Gamma(C), j \notin 
  \langle \overline{B}_{0},...,\overline{B}_{k-1} \rangle \} .\]
\hfill $\Box$

\medskip

\noindent \textbf{Proof of Proposition \ref{Inversion}}

In what follows we look at the functions as elements of 
$\hat{\mathcal{O}}_{\mathbf{C}^{2},0}$. If the local coordinates $X,Y$ are
chosen, one obtains a natural isomorphism
$\hat{\mathcal{O}}_{\mathbf{C}^{2},0} \simeq
\mathbf{C}[[X,Y]]$. If $f \in \mathbf{C}[[X]][Y]$, by definition (see
section 4), 
$B_{0} = d(f) = (f,X)$. Now, $X$ is a regular function at the
origin. We take other coordinates, $x,y \in
\hat{\mathcal{O}}_{\mathbf{C}^{2}}$ generic for the functions $f$ and
$X$. By the 
implicit function theorem, we have $C_{X} = C_{h}$, where: 
$$h(x, y) := y - \gamma(x),$$ 
with $\gamma \in \mathbf{C}[[x]]$. Take a
Newton-Puiseux parameterization of $C_{f}$: 
$$\begin{cases}
    x = T^{b_{0}} \\
    y = y(T)
\end{cases}.$$ 
Then: $B_{0} =(f,X) = (f,h) = v_{T}(h(T^{b_{0}},y(T))) = 
  v_{T}(y(T)-\gamma(T^{b_{0}}))$.

The first exponent in $y(T)$ which is not divisible by $b_{0}$ is 
$b_{1}$. So, when we vary the choice of $\gamma$, we cannot obtain
a value $v_{T}(y(T)-\gamma(T^{b_{0}}))$ greater than
$b_{1}$. The value $b_{1}$ can be obtained if the truncations
of $Y_{1}(T)$ and $\gamma(T^{b_{0}})$ coincide up to the order
$b_{1}$ (not including it). When $\gamma$ varies, we can also
obtain all the values $lb_{0}$, with $lb_{0}< b_{1}$,
i.e., with $l \leq [\frac{b_{1}}{b_{0}}].$

Once we know the degree $B_{0}$, by Proposition \ref{semigroup} we know that 
$\overline{B}_{0}=B_{0}$ and that all the numbers 
$\overline{B}_{1},...,\overline{B}_{G}$ are uniquely determined by the
minimality property from
the semigroup, which is independent of the coordinates. 
Then one can compute, in this order, the
sequences $(E_{0}, E_{1},...,E_{G})$ and $(N_{1},...,N_{G})$ and
finally obtain all the sequence $(B_{0},...,B_{G})$. 

Let us treat
successively the three cases distinguished in the statement of the
proposition. 

1) $\underline{B_{0}= b_{0}}$. 

This means that the $Y$-axis is transverse to $C_{f}$. 
  Then it is immediate that $G=g$ and 
$(\overline{B}_{0},...,\overline{B}_{G}) = 
(\overline{b}_{0},...,\overline{b}_{g})$. So:
$$ (B_{0},...,B_{G})=(b_{0},...,b_{g}).$$

2) $\underline{B_{0}= l\cdot b_{0}}$, with 
$l \in \{ 2,...,[\frac{b_{1}}{b_{0}}]\}.$

This means that the $Y$-axis is tangent to $C_{f}$ but has not maximal
contact with it (see the definition of this notion given after
Corollary \ref{dual}). Then $\overline{B}_{0}= lb_{0}$. As $b_{0}$ is
the minimal 
element of $\Gamma(C)-\{0\}$ and $b_{0}<\overline{B}_{0}$, we see
that $\overline{B}_{1} = b_{0}$. Then $E_{1}= b_{0}$ and so 
$\overline{B}_{2} = b_{1}$. Continuing like this we get: 
$$G = g+1$$
$$(\overline{B}_{0},\overline{B}_{1},\overline{B}_{2},...,\overline{B}_{G})
= (l  \overline{b}_{0},\overline{b}_{0},\overline{b}_{1},...,
\overline{b}_{g})$$
$$ (E_{0}, E_{1},E_{2},...,E_{G})= (l 
e_{0},e_{0},e_{1},
...,e_{g})$$
$$ (N_{1},N_{2},N_{3},...,N_{G}) = (l,n_{1},n_{2},...,n_{g}).$$

By proposition \ref{semigroup}, point 1), we get : $B_{k}- B_{k-1}= 
\overline{B}_{k} - N_{k-1}\overline{B}_{k-1} = \linebreak =
\overline{b}_{k-1}-N_{k-1}\overline{b}_{k-2}$, for all 
$k \in \{1,...,G\}$.

For $k=2$, $B_{2}- B_{1}=\overline{b}_{1}-l\overline{b}_{0}= 
b_{1}-lb_{0}$ which gives:
$$B_{2} =  b_{1}+(1-l)b_{0}.$$

For $k\geq 3$, $N_{k-1}=n_{k-2}$ and so:
$$ B_{k}- B_{k-1}=
\overline{b}_{k-1}-n_{k-2}\overline{b}_{k-2} = 
b_{k-1}-b_{k-2}.$$

We obtain by induction:
$$ (B_{0},...,B_{G})= (lb_{0}, b_{0}, 
b_{1}+(1-l)b_{0},..., b_{g}+(1-l)b_{0}).$$

3) $\underline{B_{0}= b_{1}}$. 

This means that the $Y$-axis has maximal contact with the branch
$C_{f}$.  
The same kind of analysis as before shows that: 
$$G=g$$ 
$$(\overline{B}_{0},\overline{B}_{1},\overline{B}_{2},...,\overline{B}_{G})
= ( \overline{b}_{1},\overline{b}_{0},\overline{b}_{2},...,
\overline{b}_{g})$$
$$ (E_{0}, E_{1},E_{2},...,E_{G})= (b_{1},e_{1},e_{2},...,
e_{g})$$
$$ (N_{1},N_{2},N_{3},...,N_{G}) =
(\frac{b_{1}}{e_{1}},n_{2},n_{3},...,n_{g})$$
$$ (B_{0},...,B_{G})= (b_{1}, b_{0}, 
b_{2}+b_{0}-b_{1},..., b_{g}+b_{0}-b_{1}).$$

\medskip

In order to deal with truncations we use Proposition \ref{Noether}. Since
we have two systems of coordinates, we note by $K^{(X,Y)}(f,\phi)$ the
coincidence exponent of $f$ and $\phi$ in the coordinates $(X,Y)$. See 
Proposition \ref{Noether} for its definition. 

Let $\phi_{k+\epsilon}$ be the $(k+\epsilon)$-semiroot of $f$ with
respect to $(X,Y)$ which is equal to the minimal polynomial of a
$(k+\epsilon)$-truncated Newton-Puiseux series of $f$ (see
Proposition \ref{existence}). Then we look at $f$ and $\phi_{k+\epsilon}$ in
the coordinates $(x,y)$ and we compute $K^{(x,y)}
(f,\phi_{k+\epsilon})$ using Proposition \ref{Noether}. This shows that some
precisely determined truncations of their Newton-Puiseux series in these
coordinates coincide. As $\phi_{k+\epsilon}$ is determined only by the 
$(k+\epsilon)$-truncation of the Newton-Puiseux series of $f$ with respect to
$(X,Y)$, the computations done in each of the three cases give the result.

We give as an example only the treatment of the second case 
$(B_{0}= lb_{0}).$

In this case $\epsilon=1$. Let us consider $k \geq 1$ and the semiroot 
$\phi_{k+1}$. We know, by Theorem \ref{Mainthm}, that $(X,\phi_{k+1})
= d(\phi_{k+1}) = \frac{N}{E_{k+1}}$ and 
$(f, \phi_{k+1}) = \overline{B}_{k+2}$. Then: $(x,\phi_{k+1})= m(\phi_{k+1}) = 
\frac{B_{1}}{B_{0}}\cdot \frac{N}{E_{k+1}} = \frac{B_{1}}{E_{k+1}} = 
\frac{b_{0}}{e_{k}}.$ So:
$$ \frac{(f, \phi_{k+1})}{(x,\phi_{k+1})} = 
 \frac{\overline{B}_{k+2}e_{k}}{b_{0}} = 
\frac{\overline{b}_{k+1}e_{k}}{b_{0}} = 
\frac{\overline{b}_{k+1}}{n_{1}\cdots n_{k}} $$
and Proposition \ref{Noether} applied in the coordinate system
$(x,y)$ gives the equality 
$K^{(x,y)}(f,\phi_{k+1}) = \frac{b_{k+1}}{n}$,
which shows that $f$ and $\phi_{k+1}$ have coinciding $k$-truncated
Newton-Puiseux series in the coordinates $(x,y)$. \hfill $\Box$

\medskip
\noindent \textbf{Proof of Proposition  \ref{exp}}

We consider expansions of the type (\ref{expansion}):
$$P = a_{0}Q^{s}+a_{1}Q^{s-1}+\cdots+a_{s}$$
with $d(a_{i})<d(Q)$ for all $i\in\{0,...,s\}$.

Let us  show  that in such an expansion, \textit{the degrees of the terms are
all different}. More precisely, we show that: 
\begin{equation} \label{inegcr}
d(a_{i}Q^{s-i}) >  d(a_{j}Q^{s-j}) \mbox{ for } i<j.
\end{equation} 

Indeed, we have:
$$d(a_{i}Q^{s-i})-d(a_{j}Q^{s-j}) = d(a_{i}) - d(a_{j}) + (j-i) d(Q)
\geq d(a_{i}) - d(a_{j}) + d(Q) > 0.$$

From this property of the degrees, one deduces that a $Q$-adic
expansion of $0 \in A[Y]$ is necessarily trivial, which in turn gives
\textit{the unicity of the expansion} for all monic $P \in A[Y]$.

Moreover, identifying the leading coefficients in both sides of
equation (\ref{expansion}), we see that $a_{0}$ is monic. 

Then we have also: $d(P) = d(a_{0}Q^{s}) = d( a_{0}) +
s d(Q)$, which implies: 
$$ \frac{d(P)}{d(Q)} = s + \frac{d(a_{0})}{d(Q)}. $$
But $0 \leq \frac{d(a_{0})}{d(Q)}< 1$, which gives the equality 
$s=[\frac{d(P)}{d(Q)}]$. Also, since $a_{0}$ is monic,
$d(Q) \mid  d(P)
\Leftrightarrow  d(a_{0})=0 \Leftrightarrow a_{0} = 1$.

Let us suppose now we are in the case when $d(Q) \mid
d(P)$. We have just
seen that in this situation $a_{0}=1$. Then:
$$ P-Q^{s} = \sum_{i=1}^{s}a_{i}Q^{s-i} $$
and, by the growth property of the degrees (\ref{inegcr}), we get:
$$ d(P-Q^{s}) \leq d(a_{2}Q^{s-2}) \Leftrightarrow a_{1}=0. $$

But $d(a_{2}Q^{s-2}) = d(a_{2})+(s-2)d(Q) < (s-1)d(Q) =
d(P)-\frac{d(P)}{d(Q)}$, and so:
$$ d(P-Q^{s}) < d(P)-\frac{d(P)}{d(Q)} \Leftrightarrow a_{1}=0. $$

By the definition of approximate roots, we see that $a_{1}=0$ if and
only if $Q = \sqrt[s]{P}$. \hfill $\Box$

\medskip

\noindent \textbf{Proof of Proposition \ref{Iterate}}

Let us take for $Q$ a monic polynomial, $d(Q)=\frac{d(P)}{p}$. The
$Q$-adic expansion
of $P$ is of the form:

$$ P = Q^{p}+a_{1}Q^{p-1}+\cdots +a_{p} $$
with $d(a_{i})<d(Q)$ for $1 \leq i \leq p$.

We consider also the $\tau_{P}(Q)$-adic expansion of $P$:
$$ P = \tau_{P}(Q)^{p}+a_{1}'\tau_{P}(Q)^{p-1}+\cdots +a_{p}'. $$

We shall prove that if $a_{1} \neq 0$, we have $d(a_{1}') <
d(a_{1})$. This will show that after iterating $\tau_{P}$ at most
$d(a_{1})+1$ times, we arrive at the situation $a_{1}=0$, in which
case $\tau_{P}(Q)=Q= \sqrt[p]{P}$. But $d(a_{1})+1 \leq d(Q) =
\frac{d(P)}{p}$, which proves the proposition. 

So, let us suppose $a_{1} \neq 0$. Then:

\begin{equation}
   \label{devpt}
      P = (Q +\frac{1}{p}a_{1})^{p} + \sum_{k=2}^{p}a_{k}Q^{p-k} - 
            \sum_{k=2}^{p}\left( \begin{array}{l}
                                       p \\ k
                                   \end{array}
                                 \right)
                                 \frac{1}{p^{k}}a_{1}^{k}Q^{p-k}.
\end{equation}

We study now the $\tau_{P}(Q)$-adic expansion of
$P-\tau_{P}(Q)^{p}$ starting from equation (\ref{devpt}). First, for $2
\leq k \leq p$, we have:
$$ d(a_{k}Q^{p-k}) < d(Q) + (p-k)d(Q) \leq (p-1)d(Q) $$
$$ d(a_{1}^{k}Q^{p-k}) < k d(Q) + (p-k)d(Q) = p\cdot d(Q). $$

But $d(\tau_{P}(Q)) =d(Q)$ and Proposition \ref{exp} shows that the
$\tau_{P}(Q)$-adic expansion of $a_{k}Q^{p-k}$ has non-zero terms of the form 
$c_{j}\tau_{P}(Q)^{j}$ with $j \leq p-2$ and the $\tau_{P}(Q)$-adic 
expansion of $a_{1}^{k}Q^{p-k}$ of the form $c_{j}\tau_{P}(Q)^{j}$
with $j \leq p-1$.

Let $c_{0}^{(k)}\tau_{P}(Q)^{p-1}$ be the term corresponding to 
$\tau_{P}(Q)^{p-1}$ in the $\tau_{P}(Q)$-adic expansion of
$a_{1}^{k}Q^{p-k}$. It is possible that $c_{0}^{(k)}=0$. Then:

$$ d(c_{0}^{(k)}\tau_{P}(Q)^{p-1}) \leq d(a_{1}^{k}Q^{p-k}) $$
and so:
$$ d(c_{0}^{(k)}) \leq k\cdot d(a_{1}) -k\cdot d(Q)+  d(Q) \leq 
                       2d(a_{1}) - d(Q) \leq d(a_{1}) -1. $$
  
But the polynomial $a_{1}'$ is a linear combination with coefficients
in $A$ of the polynomials $c_{0}^{(k)}$, which shows the announced
inequality:
$$ d(a_{1}') \leq  d(a_{1}) -1. $$

With this the proof is complete. \hfill $\Box$

\medskip

\noindent \textbf{Proof of Proposition \ref{Noether}}

As stated in section 6, this result is classical. Recent proofs are
contained in \cite{M 77}  
(for generic coordinates) and \cite{GP 91} (for arbitrary
coordinates). We give here a rather detailed proof in order to explain
the origin of the formula for $\overline{B}_{k}$ in
Proposition \ref{semigroup}.

Let $N=d(f), M=d(\phi)$. Decompose $\phi \in \mathbf{C}[[X]][Y]$ as a
product of terms of degree 1:
$$ \phi(X,Y) = \prod_{i=1}^{M}(Y-\zeta_{i}(X))$$
where the $\zeta_{i}(X)$ are all the Newton-Puiseux series of $\phi$ with
respect to $(X,Y)$.

Let $T \rightarrow (T^{N},Y(T))$ be a parameterization of $f$, obtained
from a fixed Newton-Puiseux series $\eta(X)$. As $T = X^{\frac{1}{N}}$, we
have: $Y(T) = \eta(X)$. Then, using the rules explained in section 4:
\begin{equation*}
  \begin{array}{ll}
    (f,\phi) & = v_{T}(\phi(T^{N},Y(T))) = 
                   v_{T}(\prod_{i=1}^{M}(Y(T)- \zeta_{i}(T^{N})))= \\
          & = v_{X^{\frac{1}{N}}}
                 (\prod_{i=1}^{M}( \eta(X)-\zeta_{i}(X)))= 
              N v_{X}(\prod_{i=1}^{M}( \eta(X)-\zeta_{i}(X)))= \\
          & = N \sum_{i=1}^{M}v_{X}( \eta(X)-\zeta_{i}(X))) 
  \end{array}
\end{equation*}
           
Now we look at the possible values of $v_{X}( \eta(X)-\zeta_{i}(X)))$
when $i$ varies, and for a fixed value we look how many times it is
obtained.

If $k$ is minimal such that $K(f,\phi) < \frac{B_{k+1}}{N}$, we get:

$\bullet$ the value $\frac{B_{i}}{N}$ is obtained $M \cdot 
\frac{E_{i-1}-E_{i}}{N}$ times, for $ i \in \{1,...,k\}$.

$\bullet$ the value $K(f,\phi)$ is obtained $M\cdot \frac{E_{k}}{N}$ times.

So:
$$(f,\phi)= N[\sum_{i=1}^{k}M \cdot \frac{E_{i-1}-E_{i}}{N}\cdot
\frac{B_{i}}{N} + M\cdot \frac{E_{k}}{N}\cdot K(f,\phi)]$$
which implies:
$$\frac{(f,\phi)}{M}=
\sum_{i=1}^{k}(E_{i-1}-E_{i})\frac{B_{i}}{N}+E_{k}\cdot K(f,\phi).$$
 
Now recall the formula for $\overline{B}_{k}$ given in Proposition  4:
\begin{equation} \label{Noet}
  \overline{B}_{k} = B_{k} + 
            \sum_{i=1}^{k-1}\frac{E_{i-1}-E_{i}}{E_{k-1}} B_{i}, 
\end{equation}
which gives:
$$\sum_{i=1}^{k-1}\frac{E_{i-1}-E_{i}}{B_{i}}= E_{k-1}\overline{B}_{k}-
E_{k}B_{k}.$$
We get:
$$\frac{(f,\phi)}{M}=\frac{E_{k-1}\overline{B}_{k}}{N}-
        \frac{E_{k}B_{k}}{N}+E_{k}K(f,\phi),$$
which is the desired formula.\hfill $\Box$

\medskip
 \textbf{Remark:} We had nothing to know about the relation of
$\overline{B}_{k}$ with the semigroup $\Gamma(C)$. We only needed the
fact it is given by formula (\ref{Noet}). See also the comments made
in the proof of Proposition \ref{semigroup}.

\medskip

\noindent \textbf{Proof of Proposition \ref{existence}}

Let $\eta_{k}(X)= h_{k}(X^{\frac{E_{k}}{N}})$. Then $h_{k}(T) \in
  \mathbf{C}[[T]]$ 
and:
$$ \begin{cases}
      X = T^{\frac{N}{E_{k}}} \\
      Y = h_{k}(T)
   \end{cases} $$
is a \textit{reduced} Newton-Puiseux parameterization of
$(\phi_{k}=0)$. So we have: 
$$ d( \phi_{k})=\frac{N}{E_{k}}.$$

Now, using Proposition  \ref{Noether}, since $K(f, \phi_{k}) =
\frac{B_{k+1}}{N}$, we have: 
$$ \frac{(f,\phi_{k})}{d(
  \phi_{k})}=\frac{\overline{B}_{k+1}}{N_{1}\cdots N_{k}} \Rightarrow 
    (f,\phi_{k}) = \frac{N}{E_{k}}\cdot 
        \frac{\overline{B}_{k+1}}{N_{1}\cdots N_{k}}=
        \overline{B}_{k+1}. $$

We have obtained: $d( \phi_{k})=\frac{N}{E_{k}}$ and 
  $(f,\phi_{k}) = \overline{B}_{k+1}$, which shows that $\phi_{k}$ is
  a $k$-semiroot. \hfill $\Box$

\medskip

\noindent \textbf{Proof of Proposition \ref{indep}}

We prove the proposition by induction on $k$. 

For $k=0$, we have $\phi \in \mathbf{C}[[X]]$, and so
$\phi(X)=X^{M}u(X)$, where 
$u(0)\neq 0$, so: $(f,\phi) = M\cdot (f,x) = M \cdot d(f) \in \langle
B_{0} \rangle = \langle\overline{B}_{0}\rangle$. 

Suppose now the proposition is true for $k \in \{0,...,G-1\}$. We 
prove it for $k+1$. 

Consider $\phi \in \mathbf{C}[[X]][Y], d(\phi) < \frac{N}{E_{k+1}}$ and take
a $k$-semiroot $q_{k}$ of $f$, which exists by Proposition \ref{existence}. 
Make the $q_{k}$-adic expansion of $\phi$:
$$ \phi = a_{0} q_{k}^{s} + a_{1} q_{k}^{s-1} +\cdots + a_{s}. $$

We prove that the intersection numbers $(f,a_{i} q_{k}^{s-i})$ are all
distinct. Suppose by contradiction that $0 \leq j < i \leq s$ and 
$(f,a_{i} q_{k}^{s-i})=(f,a_{j} q_{k}^{s-j})$. Then $(i-j)(f, q_{k}) = 
\linebreak =
(f,a_{i}) - (f,a_{j}) \in \langle\overline{B}_{0},...,
\overline{B}_{k}\rangle$, by the induction hypothesis. So: 
$E_{k} \mid  (i-j)\overline{B}_{k+1}$. But $E_{k+1} = \mbox{gcd}(E_{k},
\overline{B}_{k+1})$, and so we obtain: $\frac{E_{k}}{E_{k+1}}\mid
(i-j)$. Now, $\frac{E_{k}}{E_{k+1}}=N_{k+1}$ and $i-j \leq s= 
[\frac{d(\phi)}{d(q_{k})}]$. As $ \frac{d(\phi)}{d(q_{k})} =
\frac{E_{k}}{N}\cdot d(\phi) < \frac{E_{k}}{N}\cdot \frac{N}{E_{k+1}} =
N_{k+1}$, we see that $s < N_{k+1}$, which gives a contradiction.

This shows that the numbers $(f,a_{i} q_{k}^{s-i})$ are all distinct
and so:
$$ (f,\phi) = \mbox{min}_{i}\{(f,a_{i} q_{k}^{s-i})\} = 
    \mbox{min}_{i} \{(f,a_{i})+ (s-i)(f, q_{k})\}. $$
But $(f,a_{i}) \in \langle\overline{B}_{0},...,
\overline{B}_{k}\rangle$ by the induction hypothesis and $(f,\phi) = 
\overline{B}_{k+1}$, so: $(f,\phi) \in \langle\overline{B}_{0},...,
\overline{B}_{k+1}\rangle$. With this, the step of induction is
completed. \hfill $\Box$

\medskip

\noindent \textbf{Proof of Proposition \ref{stab}}

If $\phi$ is a $k$-semiroot and $\psi$ a $(k-1)$-semiroot, then 
$d(\phi)= \frac{N}{E_{k}}$ and $d(\psi)= \frac{N}{E_{k-1}}$. So the
$\psi$-expansion of $\phi$ is of the form:
\begin{equation} \label{rootroot}
   \phi = \psi^{N_{k}}+ a_{1}\psi^{N_{k}-1}+\cdots +a_{N_{k}}.
\end{equation}

We have $\tau_{\phi}(\psi)= \psi+ \frac{1}{N_{k}}a_{1}$.

We are going to show that: 
$$(f,\psi) < (f,a_{1}).$$ 

This will give $(f, \tau_{\phi}(\psi))= (f,\psi)=\overline{B}_{k}.$ But
$d(a_{1}) < d(\psi)$ and so $d(\tau_{\phi}(\psi))=
\linebreak =d(\psi)$, which shows that
$\tau_{\phi}(\psi)$ is also a $(k-1)$-semiroot.

Exactly as in the proof of Proposition \ref{indep}, we have that the
intersection numbers $(f,
a_{i}\psi^{N_{k}-i})$ are all distinct, for $i \in \{1,...,N_{k}\}.$
Using equation (\ref{rootroot}) we deduce:

$$(f,\phi-\psi^{N_{k}})= \mbox{min}_{1 \leq i \leq N_{k}}
\{(f,a_{i}\psi^{N_{k}-i})\} 
\leq (f,a_{1}\psi^{N_{k}-1}).$$ 

But, by Proposition \ref{semigroup}, $(f,\phi)=\overline{B}_{k+1}> 
N_{k}\overline{B}_{k}= (f,\psi^{N_{k}})$ and so: \linebreak 
$(f,\phi-\psi^{N_{k}})=(f,\psi^{N_{k}})$. We obtain:
$$  (f,a_{1})+( N_{k}-1)(f,\psi) \geq N_{k}(f,\psi) $$
which gives:
$$(f,\psi)\leq (f,a_{1}).$$

On the other hand, $d(a_{1})<d(\psi) = \frac{N}{E_{k-1}}$, and
Proposition \ref{indep} 
shows that $(f,a_{1}) \in \langle\overline{B}_{0},...,
\overline{B}_{k-1}\rangle$. But $(f,\psi) = \overline{B}_{k}\notin 
\langle\overline{B}_{0},...,\overline{B}_{k-1}\rangle$, which shows
that we cannot have the equality $(f,\psi)=(f,a_{1})$.

Thus, we have proven the inequality $(f,\psi) <(f,a_{1})$ and with it
the proposition.  \hfill $\Box$

\section{The Approximate Roots \\and the Embedding Line Theorem}

We present the ideas of the proofs of the epimorphism theorem and of
the embedding line theorem as they are given in \cite{A 77}. 

It is in order to do these proofs that is developed in \cite{A 77} the
theory of Newton-Puiseux 
parameterizations and of local semigroups for elements of
$\mathbf{C}((X))[Y]$, the {\it meromorphic curves}. This
framework is more 
general than the one presented before, which concerned elements of 
$\mathbf{C}[[X]][Y]$, the {\it entire curves}. We have chosen to give  
before all the proofs for entire curves, first because they are in
general used for the
local study of plane curves and second in order to
point out in this final section the differences between the two
theories. A third type of curves, the {\it purely
  meromorphic} ones, will prove to be of the first importance. 

\medskip

\noindent \textbf{Proof of the Epimorphism Theorem}

We consider an \textit{epimorphism} $\sigma: \mathbf{C}[X,Y] \rightarrow
\mathbf{C}[T]$ and we note:
$$ P(T):= \sigma(X) , Q(T):= \sigma(Y), $$
$$ N:=d_{T}(P), M:=d_{T}(Q). $$

We suppose that both degrees are non zero.
The ideal $\mbox{ker}(\sigma)$ is of height one in $\mathbf{C}[X,Y]$,
so it is generated by one element. A privileged generator is given by:
$$ F(X,Y)=\mbox{Res}_{T}(P(T)-X,Q(T)-Y). $$

Here $\mbox{Res}_{T}$ denotes the \textit{resultant} of the two polynomials,
seen as polynomials in the variable $T$.

From the determinant formula for the resultant, we obtain:
$$d_{X}(F)=M, d_{Y}(F)=N $$
and that $F$ is monic if we see it as a polynomial in $X$ or in $Y$.

Let us consider the set:
$$\Gamma(F): = \{d_{T}(G(P(T),Q(T))), G \in \mathbf{C}[X,Y] - (F)
\}.$$

The set $\Gamma(F)$ is a sub-semigroup of $(\mathbf{N},+)$. \textit{The
morphism $\sigma$ is an epimorphism if and only if} $T
\in \mbox{im}(\sigma)$, \textit{which is equivalent to $1 \in \Gamma(F)$, or 
$\Gamma(F)=\mathbf{N}$}.

Make now the change of variables: $x=X^{-1}, y=Y$. Take:
$$f(x,y) := F(x^{-1},y) \in \mathbf{C}[x^{-1}][y].$$

The polynomial $f$ is monic in $y$, of degree $d(f)=N$. By
definition, the elements of $\mathbf{C}[x^{-1}][y]$ are called {\it
  purely meromorphic curves} (notation 
of \cite{A 77}). As we have the embedding of rings $\mathbf{C}[x^{-1}]
\hookrightarrow \mathbf{C}((x))$, we can also look at $f$ as being a
\textit{meromorphic curve}, i.e. an element of $ \mathbf{C}((x))[y]$.
The theory of Newton-Puiseux expansions can be generalized to elements of 
$ \mathbf{C}((x))[y]$, and so $f$ has associated Newton-Puiseux series
$\eta(x) \in \mathbf{C}((x^{\frac{1}{N}}))$ and Newton-Puiseux
parameterizations of the form: $x=\tau^{N}, y=y(\tau) \in
\mathbf{C}((\tau)).$ It is important here that the exponent of $\tau$
in $x(\tau)$ is taken positive (see below).

From such a primitive Newton-Puiseux parameterization (see the definition in
section 2), one can obtain a  characteristic sequence of
integers $(B_{0},...,B_{G})$, where we put $B_{0}=-N$ and the other
$B_{i}$'s are defined recursively as in the case of $ \mathbf{C}[[x]][y]$,
treated before. At the same time we define the sequence of greatest
common divisors $(E_{0},...,E_{G})$, which are elements of
$\mathbf{N}^*$, and the sequences $(N_{1},...,N_{G}),
(\overline{B}_{0},...,\overline{B}_{G})$, as in section 2. Notice
that $(B_{0},...,B_{G})$ is again a strictly increasing sequence, but
not necessarily $(\overline{B}_{0},...,\overline{B}_{G})$.

If $\phi \in \mathbf{C}((x))[y]$, $f \not\vert \phi$, we define:
$$ (f,\phi) := v_{x}(\mbox{Res}_{y}(f,\phi)).$$

This construction extends the definition of the intersection number from 
$\mathbf{C}[[x]][y]$ to $\mathbf{C}((x))[y]$. It is again true with
this definition that:
$$(f,\phi)=v_{\tau}(\phi(\tau^{N},y(\tau))),$$
if $\tau \rightarrow (\tau^{N},y(\tau))$ is a Newton-Puiseux parameterization
of $f$ (we understand here why it is important to take $x=\tau^{N}$
and not $x=\tau^{-N}$).

We define now:
$$\Gamma_{\mathbf{C}[x^{-1}]}(f) := 
    \{(f,\phi),\phi \in \mathbf{C}[x^{-1}][y], f\not\vert \phi \}.$$

The set $ \Gamma_{\mathbf{C}[x^{-1}]}(f)$ is a sub-semigroup of
$(\mathbf{Z},+)$. In fact we can say more. Indeed, if
$\Phi(X,Y)=\phi(X^{-1},Y) \in \mathbf{C}[X,Y]$, we have:
$$d_{T}(\Phi(P(T),Q(T)))= -(f,\phi),$$
which shows that:
$$\Gamma_{\mathbf{C}[x^{-1}]}(f) = -\Gamma(F).$$ 

We see in particular that the semigroup
$\Gamma_{\mathbf{C}[x^{-1}]}(f)$ consists only of negative numbers.

As $\sigma$ is an epimorphism, we get:
$$\Gamma_{\mathbf{C}[x^{-1}]}(f) = \mathbf{Z}_{-}.$$

\textbf{Remark:} If we consider intersections with elements of 
$\mathbf{C}((x))[y]$, we can define a second semigroup 
$\Gamma_{\mathbf{C}((x))}(f)$. We have obviously the inclusion 
$\Gamma_{\mathbf{C}[x^{-1}]} \subset \Gamma_{\mathbf{C}((x))}$,
but in general 
this is not an equality. 

Take for example $f= y^{2}-x^{-1}$. A Newton-Puiseux parameterization of $f$
is $\tau \rightarrow (\tau^{2}, \tau^{-1})$. Take $\phi =
y^{2}-(x^{-1}-x) \in \mathbf{C}((x))[y] -
\mathbf{C}[x^{-1}][y]$. Compute their intersection number:
$$(f,\phi) = v_{\tau}(\phi(\tau^{2}, \tau^{-1})) = 2 \notin \mathbf{Z}_{-} = 
\Gamma_{\mathbf{C}[x^{-1}]}(f).$$

\medskip

Suppose now by contradiction that we are in a case where neither
$N\mid M$ nor $M \mid N$. This implies easily that $B_{1}=-M$. Indeed, 
$v_{\tau}(y)=(f,y)= \linebreak =- d_{T}(Y(P(T),Q(T)))= -
d_{T}(Q(T))=-M.$ Since $N\not\vert M$ we deduce by the definition
of $B_{1}$ that $B_{1}=v_{\tau}(y)=-M.$

Since $\Gamma_{\mathbf{C}[x^{-1}]}(f)= \mathbf{Z}_{-}$, we get in
particular $-E_{1} \in \Gamma_{\mathbf{C}[x^{-1}]}(f)$.

The contradiction is got in \cite{A 77} from the properties:
$$B_{0}=-N,\: B_{1}=-M,\: -E_{1} \in \Gamma_{\mathbf{C}[x^{-1}]}(f).$$

\medskip

\textit{Here is the place in the proof where the approximate roots make their
appearance.} As in the case of $\mathbf{C}[[X]]$, from the sequences
$(B_{0},...,B_{G})$ and $(E_{0},...,E_{G})$ one can define inductively
a sequence $(\overline{B}_{0},...,\overline{B}_{G})$ by the relations
given in Proposition \ref{semigroup}.

They are elements of $\Gamma_{\mathbf{C}((x))}(f)$, as they can be
obtained by intersecting $f$ with arbitrary semiroots of $f$, for
example the ones got by truncating a Newton-Puiseux series of $f$
(Proposition \ref{existence}
generalizes to this context).

But, more important, $(\overline{B}_{0},...,\overline{B}_{G})$ {\it are
elements of} $\Gamma_{\mathbf{C}[x^{-1}]}(f)$. Indeed, $f_{k}=
\linebreak =
\sqrt[E_{k}]{f} \in \mathbf{C}[x^{-1}][Y]$. Theorem \ref{Mainthm}
generalizes to 
this context and so: $(f,f_{k})= \linebreak =\overline{B}_{k+1}$, for $k \in
\{0,...,G\}$.

What is again true is that  $(\overline{B}_{0},...,\overline{B}_{G})$
form {\it a system of generators} of
$\Gamma_{\mathbf{C}[x^{-1}]}(f)$. As $\Gamma_{\mathbf{C}[x^{-1}]}(f)$
is composed of negative numbers, we cannot speak any more about a
minimal system of generators, as in Proposition \ref{semigroup}. 
What remains true is that
they are {\it a strict system of generators} (see \cite{A 77}) in the
following sense:

\begin{prop} \label{strict}
 Every element $\gamma$ of $\Gamma_{\mathbf{C}[x^{-1}]}(f)$ can be
 expressed in a unique way as a sum:
$$\gamma= i_{-1}\overline{B}_{0}+\cdots +i_{G-1}\overline{B}_{G}$$
where $i_{-1} \in \mathbf{N}$ and $0 \leq i_{k} <N_{k+1}$ for $k \in
\{1,...,G-1\}.$
\end{prop}

\medskip

To get this proposition, one proves first an analog of Corollary
\ref{expansions}, obtained by replacing $\mathbf{C}[[X]]$ by
$\mathbf{C}[x^{-1}]$. The proof follows the same path.

Now write the property $-E_{1} \in \Gamma_{\mathbf{C}[x^{-1}]}(f)$ in
terms of this strict sequence of generators:
$$ -E_{1} = i_{-1}\overline{B}_{0}+\cdots +i_{G-1}\overline{B}_{G}.$$

Take $p := \mbox{max}\{k \in \{0,...,G\}, i_{k-1}\neq 0\}.$ So:
$$
 E_{1} = i_{-1}\mid \overline{B}_{0}\mid +\cdots +i_{p-1}\mid
\overline{B}_{p}\mid,$$
with $i_{p-1}\neq 0$. 

If $p\geq 2$, we get: $ E_{p-1} \mid (E_{1}-i_{-1}\mid
\overline{B}_{0}\mid - \cdots -i_{p-2}\mid
\overline{B}_{p-1}\mid)$ and so: $E_{p-1} \mid (i_{p-1}\mid
\overline{B}_{p}\mid).$ Since $E_{p}=\mbox{gcd}(E_{p-1},\mid
\overline{B}_{p}\mid)$, we get $N_{p}=\frac{E_{p-1}}{E_{p}} \mid
i_{p-1}$, which contradicts the inequality $0 < i_{p-1}<N_{p}$.

So we obtain $p\leq 1$ and:
$$ E_{1} = i_{-1}\mid \overline{B}_{0}\mid +i_{0}\mid
\overline{B}_{1}\mid.$$

This implies: $1 = i_{-1}\frac{\mid\overline{B}_{0}\mid}{E_{1}}+
i_{0}\frac{\mid\overline{B}_{1}\mid}{E_{1}}$, which shows that 
$\frac{\mid\overline{B}_{0}\mid}{E_{1}}=1$ or 
$\frac{\mid\overline{B}_{1}\mid}{E_{1}}=1$. But, by the recursive
relations giving the $\overline{B}_{i}$'s, $\overline{B}_{0}=B_{0}=-N$
and $\overline{B}_{1}=B_{1}=-M$, so:
$$ \frac{N}{(M,N)}=1 \mbox{ or }\frac{M}{(M,N)}=1.$$

We get: $M\mid N$ or $N\mid M$, which contradicts our hypothesis. The
theorem is proved. \hfill $\Box$

\medskip

 \textbf{Remark:} One can also give a  proof without using
  contradiction. In this case one cannot suppose from the beginning
  that $N \not| M$, and so it is not necessarily true that
  $B_{1}=-M$. As one cannot hope to express in this case $B_{1}$ in
  terms of $N$ and $M$, the preceding proof appears to get in
  trouble. This can be arranged if one modifies the definition of the
  characteristic sequence, taking for $B_{1}$ the minimal
  exponent appearing in $y(\tau)$, without imposing that it should not be
  divisible by $N$. This is the definition of characteristic sequence
  taken in a majority of Abhyankar's writings on curves, in particular
  \cite{A 77}, where the preceding proof is given with this modified
  definition.
 
\medskip

\noindent \textbf{Proof of the Embedding Line Theorem}

If the epimorphism $\sigma: \mathbf{C}[X,Y] \rightarrow
\mathbf{C}[T]$ is given by $X=P(T), Y= Q(T)$, put $N := d_{T}(P), 
M:=d_{T}(Q)$ and write:
$$\begin{array}{l}
P(T) =\alpha_{0}T^{N}+\alpha_{1}T^{N-1}+\cdots+\alpha_{N},\\
Q(T) =\beta_{0}T^{M}+\beta_{1}T^{M-1}+\cdots+\beta_{M}.
\end{array}$$
(we consider here that $d_{T}(0)=0$).

Suppose one of the degrees $M,N$ is zero, for example $M=0$. Then:
$Q(T)=\linebreak =\beta_{0} \in \mathbf{C}$.

For all $G \in \mathbf{C}[X,Y]$, $d_{T}(G(P(T),Q(T)) \in
N\mathbf{N}.$ If $\sigma$ is an epimorphism, there exists such a $G$
with $d_{T}(G(P(T),Q(T))=1$, and this implies $N=1$. So:
$$\begin{cases}
     P(T)=\alpha_{0}T+\alpha_{1}, \; \alpha_{0}\neq 0 \\
     Q(T)= \beta_{0}
\end{cases}.$$

Consider the isomorphism of $C$-algebras $\sigma_{1}:\mathbf{C}[U,V]
\rightarrow \mathbf{C}[X,Y]$, given by:
$$\begin{cases}
     U=\frac{1}{\alpha_{0}}X-\frac{\alpha_{1}}{\alpha_{0}} \\
     V = Y-\beta_{0}
\end{cases}.$$

Then $\sigma \circ \sigma_{1}: \mathbf{C}[U,V]
\rightarrow \mathbf{C}[T]$ is given by:
$$\begin{cases}
     U= T \\
     V = 0
\end{cases}$$
and the theorem is proved in this case.

\medskip

Suppose now that $M \geq 1, N \geq 1$. \textit{By the epimorphism theorem},
$M\mid N$ or $N \mid M$. Suppose for example that $M \mid N$. Consider
the isomorphism of $\mathbf{C}$-algebras $\sigma_{1}:\mathbf{C}[U,V]
\rightarrow \mathbf{C}[X,Y]$ given by:
$$\begin{cases}
     U= X - \alpha_{0}\beta_{0}^{-\frac{N}{M}}Y^{ \frac{N}{M}}\\
     V = Y
\end{cases}.$$

Then $\sigma \circ \sigma_{1}$ is given by:
$$\begin{cases}
     U= P(T)- \alpha_{0}\beta_{0}^{-\frac{N}{M}}Q(T)^{ \frac{N}{M}}\\
     V = Q(T)
\end{cases}.$$

We have: $d_{T}(P(T)- \alpha_{0}\beta_{0}^{-\frac{N}{M}}Q(T)^{
  \frac{N}{M}}) < d_{T}(P(T))$ and so in the new coordinates
  $(U,V)$, the sum of the degrees of the polynomials giving the
  embedding of the line in the plane is strictly less than in the
  coordinates $(X,Y)$.

Repeating this process a finite number of times, we see that we arrive
at the situation where one of the polynomials is a constant, the case 
first treated. This proves the theorem. \hfill $\Box$

\end{document}